\documentclass[letterpaper,11pt]{amsart}
\usepackage{amsmath,amsthm,amssymb,mathtools,amsfonts,mathrsfs}
\usepackage{xcolor}
\usepackage[breaklinks,colorlinks,citecolor=teal,linkcolor=teal,
  urlcolor=teal,pagebackref,hyperindex]{hyperref}
\usepackage[alphabetic]{amsrefs}
\usepackage{tikz-cd}
\usepackage{graphicx}
\newtheorem{theorem}{Theorem}[section]
\newtheorem{proposition}[theorem]{Proposition}
\newtheorem{corollary}[theorem]{Corollary}
\newtheorem{lemma}[theorem]{Lemma}
\theoremstyle{definition}
\newtheorem{remark}[theorem]{Remark}

\newcommand{\PP}{\mathbb P}
\renewcommand{\AA}{\mathbb A}
\newcommand{\QQ}{\mathbb Q}
\newcommand{\OO}{\mathcal O}

\newcommand{\cR}{\mathcal R}

\newcommand{\cV}{\mathcal V}
\newcommand{\scrC}{\mathscr C}
\newcommand{\vir}{\mathrm{vir}}
\newcommand{\ev}{\operatorname{ev}}
\newcommand{\pr}{\operatorname{pr}}
\newcommand{\Bl}{\operatorname{Bl}}

\title[Gromov--Witten theory of smooth pairs via projective bundles]
{Gromov--Witten theory of smooth pairs via projective bundles}
\author{Yu Wang and Fenglong You}
\address{Department of Mathematics \\ Louisiana State University \\303 Lockett Hall \\ Baton Rouge \\ LA \\ United States of America \\ 70803}
\email{yuwang@lsu.edu}

\address{School of Mathmatical Sciences \\ University of Nottingham, \\University Park \\Nottingham, NG7 2RD \\United Kingdom.}
\email{fenglong.you@nottingham.ac.uk}

\begin{document}
\begin{abstract}
We express genus zero logarithmic Gromov-Witten invariants of a smooth
pair $(X, D)$ in terms of logarithmic Gromov--Witten invariants of the corresponding projective bundle \(P=\mathbb P_X(\mathcal O_X(-D)\oplus\mathcal O_X)\) and 
one of the  relative markings becomes an interior marking. This may be considered as a logarithmic version of the main result of Wang--You \cite{WY} without a nefness assumption on $D$. 
As a corollary, we also have a log--orbifold correspondence for these invariants of the projective bundle \(P\) when $D$ is nef.
\end{abstract}
\maketitle
\tableofcontents
\clearpage

\section{Introduction}\label{sec:intro}

Let $X$ be a smooth projective variety and $D$ be a smooth divisor of $X$. The local-relative correspondence, or the local-log principle,  was first observed by \cite{Takahashi} which relates the genus zero maximal contact relative Gromov--Witten invariants of  $(\mathbb P^2,E)$ and the local Gromov--Witten invariants of $K_{\mathbb P^2}$. This example was proved by \cite{Gathmann03}. For a general smooth nef pair, it was proved by \cite{vGGR}. Some generalization in the maximal contact and nef case includes \cite{BBvG24}, \cite{BBvG21}, \cite{TY20b}, \cite{BFGW} and \cite{BNTY}.

In \cite{WY}*{Theorem 1.3}, we study Gromov--Witten invariants beyond maximal contacts where relative invariants of a smooth nef
pair $(X,D)$ are related to orbifold invariants of a projective
bundle over $X$. On the level of invariants, \cite{WY}*{Theorem 1.3} can be stated as the following identity:
\begin{equation}\label{iden-rel-orb}
\begin{split}
&\langle [\iota^*\phi_0]_{k_0},[\phi_1]_{k_1},\ldots,[\phi_n]_{k_n}\rangle_{0,n+1,\beta}^{(X,D)}\\
=&(-1)^{k_0}\langle [\phi_0\cup(X_\infty-\pi^*D)]_{(0,0)}, [\phi_1]_{(k_1,k_1)},\ldots,[\phi_n]_{(k_n,k_n)} \rangle_{0,n+1,\beta+(\sum_{i=1}^n k_i)f}^{(P,X_{\infty}+X_{\sigma})},
\end{split}
\end{equation}
where $$\pi:P:=\mathbb P(\mathcal O_X(-D)\oplus \mathcal O_X)\rightarrow X$$ is a $\mathbb P^1$-bundle over $X$ where $X_{\infty}$ is the infinity divisor of $P$ and $X_\sigma$ is the graph of the section $\sigma$ of $\mathcal O_X(D)$ that defines $D$. The left-hand side of (\ref{iden-rel-orb}) represents relative invariants of $(X,D)$ and the right-hand side of (\ref{iden-rel-orb}) are orbifold invariants of $(P,X_{\infty}+X_{\sigma})$ defined in \cite{TY20c}. Set $\sum_{i=0}^n k_i=D\cdot \beta$ for $k_i\in \mathbb Z_{>0}$, $\phi_0\in H^*(X)$ and $\phi_i\in H^*(D)$, for $i=1,\ldots,n$, as well as $\iota^*\phi_0\in H^*(D)$ which is the pullback of $\phi_0$ via the inclusion map $\iota:D\hookrightarrow X$; $f$ is the class of a fiber of the bundle $P\rightarrow X$.

As the invariants on the right-hand side of (\ref{iden-rel-orb}) are projective bundle invariants instead of local invariants, it is probably better to call Identity (\ref{iden-rel-orb}) a \emph{projective bundle-relative correspondence} or a \emph{projective-relative correspondence} instead of a local-relative correspondence. 


Some natural questions arise from the projective-relative correspondence. The first
question is whether the orbifold invariants on the right-hand side of (\ref{iden-rel-orb}) can be replaced by logarithmic invariants.
Note that, because $D$ is smooth, the relative invariants on the left-hand side of (\ref{iden-rel-orb}) are also logarithmic invariants, see \cite{AMW}. The second question is whether the nefness assumption is essential for the correspondence. As we no longer work with local invariants, the nefness of the divisor is not required to define the invariants of the projective bundle. The nefness was used in the proof of the main result of \cite{WY}. 

These two questions represent two important extensions of this framework. By working with log invariants, it laid the groundwork for the possible generalization to log invariants of snc pairs which are usually difficult to compute. Relative/log invariants with a non-nef divisor also appear naturally in degeneration and mirror symmetry. In particular, they appear in the mirror symmetry conjecture of Doran--Harder--Thompson \cite{DHT} for Tyurin degenerations. For example, the classical example of the degeneration of the quintic threefold $Q_5$ to $Q_4\cup_{K3} \operatorname{Bl} \mathbb P^3$, where $Q_4$ is a quartic threefold, $\operatorname{Bl} \mathbb P^3$ is a blow-up of $\mathbb P^3$ and they intersect along the quartic $K3$ surface. The pair $(\operatorname{Bl} \mathbb P^3, K3)$ is a non-nef pair. These relative invariants are also difficult to compute. See, for example, some recent work \cite{BGL} and \cite{You25b}. Being able to compute relative invariants of non-nef pairs allows us to compute absolute invariants through degeneration, and it is also essential in proving the Doran--Harder--Thompson conjecture.

In this paper, we prove the logarithmic version of the identity without the nefness condition assumed.
The pair on the right (\ref{iden-rel-orb}) has two boundary components, so the
comparison theorem for a smooth divisor cannot carry over directly in this scenario. We first blow up the intersection of these two
components. For the curve class, we argue that the resulting
maps avoid the strict transform, and this gives a smooth pair
$(V,E)$ for which logarithmic and expanded-pair virtual
classes introduced in \cite{Ran} have the same pushforward. We then use a degeneration formula to compare $(V,E)$ with $(X,D)$.
Thus, the new challenge comes down to analysis of what kind of rigid tropical maps contribute non-trivially to the degeneration formula, which yields the comparison result of logarithmic invariants between two smooth pairs. The preliminary logarithmic comparison
connects that calculation to the SNC target in the theorem.

\begin{remark}
    We remark that the proof of the correspondence in this paper does not rely on the proof in \cite{WY}. We work purely on the logarithmic Gromov--Witten theory without any part of the orbifold Gromov--Witten theory appearing. \cite{WY} is a generalization of the local-orbifold correspondence \cite{vGGR}, \cite{TY20b} beyond maximal contacts. The result here is a generalization of the local-log principle \cite{vGGR} beyond maximal contacts. Our result also makes it possible to explore a more general correspondence theorem for logarithmic and orbifold Gromov-Witten invariants.
\end{remark}

\subsection{The basic setup}
We work over $\mathbb C$ in genus zero, and use rational
Chow groups throughout. Curve classes are fixed in integral
homology, and rational coefficients are used only for Chow classes.

Let $X$ be a smooth projective variety with $D\subset X$ a smooth divisor, and fix an effective curve class
$\beta\in H_2(X,\mathbb Z)$ and positive
integers
\[
\vec k=(k_0,k_1,\ldots,k_n),\qquad
d=D\cdot\beta=\sum_{i=0}^n k_i>0.
\]
There are $n+1$ labeled relative markings $p_0,\ldots,p_n$ and
$n_0\geq0$ additional labeled interior markings
$p_{n+1},\ldots,p_{n+n_0}$. Put
\begin{equation}\label{eq:K}
K=\sum_{i=1}^n k_i=d-k_0.
\end{equation}
In particular, $n$ is allowed to be zero.

Consider the projective bundle
\[
\pi:P:=\PP_X\bigl(\OO_X(-D)\oplus\OO_X\bigr)\longrightarrow X.
\]
Let $\sigma\in H^0(X,\OO_X(D))$ be the section defining $D$ and $X_\sigma$ be the graph of it as introduced in \cite{WY}, that is, $X_\sigma$ is defined by 
$$(\sigma,1):\OO_X(-D)\oplus\OO_X\to\OO_X.$$
We use the quotient convention for projectivization and let $X_0$ be the zero section whose normal bundle inside $P$ is 
$\OO_X(-D)$ and $X_\infty$ be the divisor at infinity. 
The two positive sections $X_\infty,X_\sigma$ meet transversely
along a copy of $D$ and are disjoint from $X_0$. Their classes satisfy
\begin{equation}\label{eq:divisors}
[X_\sigma]=[X_\infty],\qquad
[X_0]=[X_\infty]-\pi^*[D].
\end{equation}
As in Wang--You \cite{WY}, we identify $\beta$ with its pushforward to $P$
along the section $X\simeq X_0\hookrightarrow P$.
Let $f\in H_2(P,\mathbb Z)$ be the fiber class and set
\begin{equation}\label{eq:gamma}
\gamma=\beta+Kf.
\end{equation}
Regarding $X$ and $P$ as log schemes with their respective associated divisorial logarithmic structures given by the corresponding divisors, let $\mathscr M_X$ be the moduli stack of basic stable logarithmic maps to $(X,D)$
of class $\beta$, with contact profile $\vec k$ and the $n_0$ interior
markings, and let $\mathscr M_P$ be the moduli stack for $(P,X_\infty+X_\sigma)$, of class
$\gamma$, with contact vectors
\begin{equation}\label{eq:contacts}
(0,0)\text{ at }p_0,\qquad
(k_i,k_i)\text{ at }p_i\ (1\leq i\leq n),
\end{equation}
and zero contacts at the additional interior markings where $p_0$ turns into an interior marking from a relative marking.

Write $\mathscr{\overline{M}}_{0,r}(Z,\eta)$ for the usual moduli stack of genus
zero stable maps to a smooth projective variety $Z$, with $r$ labeled markings and curve class
$\eta$. Both moduli stacks $\mathscr M_X$ and $\mathscr M_P$ possess proper maps to
\begin{equation}\label{eq:common}
\scrC=\mathscr{\overline{M}}_{0,n+n_0+1}(X,\beta)\mathbin{\times}_{X^n}D^n.
\end{equation}
Here, the component $D^n$ records the evaluations of stable logarithmic maps at $p_1,\ldots,p_n$.
The maps after pushforward to $\scrC$ forget the logarithmic structure and contact orders, and then stabilize to become 
maps to $X$ with all markings retained. Denote these two natural maps by
$\pi_{\log}$ and $u_{\log}$ respectively.
We use $\ev_i$ for evaluation at $p_i$ on the unexpanded target,
and $\ev$ for evaluation at a gluing marking. 

\begin{theorem}[Projective-log correspondence]\label{thm:main}
With the above notation, there is an identity on the level of virtual classes
\begin{equation}\label{eq:main}
(\pi_{\log})_*[\mathscr M_X]^\vir
=
(-1)^{k_0}(u_{\log})_*
\left(\ev_0^*[X_0]\cap[\mathscr M_P]^\vir\right)
\quad\text{in }A_*(\scrC)_{\QQ}.
\end{equation}
where both sides use the virtual classes of basic stable logarithmic
maps. 
\end{theorem}

Let us explain the curve class on the right-hand side of (\ref{eq:main}).
The markings $p_1,\ldots,p_n$ account for $K$ intersections
with each of $X_\infty$ and $X_\sigma$. The class $\gamma$
was chosen so that
\[
X_\infty\cdot\gamma=X_\sigma\cdot\gamma=K,\qquad
X_0\cdot\gamma=K-d=-k_0.
\]
Thus, the order $k_0$ still enters the degree, although $p_0$
has become an interior marking. For example, a curve of class $\beta$ having contact profile $(2,3)$ with $D$ corresponds to a curve of class $\gamma=\beta+3f$ where only the second relative marking is retained with a 
contact profile $(3,3)$ with $X_\infty$ and $X_\sigma$.

\begin{remark}[Insertions and descendants]\label{rem:insertions}
Since \eqref{eq:main} is an identity on $\scrC$, the insertions
at $p_1,\ldots,p_n$ may be arbitrary boundary classes.
The insertion at $p_0$ is an ambient class: on the left-hand side of (\ref{eq:main}), we
restrict a class from $X$ to $D$, and on the right-hand side of (\ref{eq:main}) we pull it
back to $P$ and multiply by $[X_0]$.
We may also apply operational classes on $\scrC$.
In particular, descendants pulled back from the common
target space $\mathscr{\overline{M}}_{0,n+n_0+1}(X,\beta)$ are allowed.
These are the descendant classes used in this paper.
\end{remark}

\begin{corollary}[log-orbifold correspondence]
\label{cor:orbifold}
Let $\mathscr M_P^{\mathrm{orb}}$ be the
moduli stack of orbifold stable maps to $(P, X_\infty+X_\sigma)$ and let  $u_{\mathrm{orb}}$ be its pushforward map to $\mathscr C$ used in Wang-You
\cite{WY}*{Theorem 1.3} with the prescribed profile shown in \eqref{eq:gamma}-\eqref{eq:contacts}. Then
\begin{equation}\label{eq:comparison}
\begin{split}
(u_{\log})_*\left(\ev_0^*[X_0]\cap[\mathscr M_P]^\vir\right)
={}&(u_{\mathrm{orb}})_*
\left(\ev_0^*[X_0]\cap[\mathscr M_P^{\mathrm{orb}}]^\vir\right)
\end{split}
\end{equation}
in $A_*(\scrC)_{\QQ}$.
\end{corollary}
\begin{proof}
For the smooth pair $(X,D)$, the logarithmic and relative virtual
classes have the same pushforward by the comparison theorem
\cite{AMW}. Theorem~\ref{thm:main} and \cite{WY}*{Theorem 1.3}
identify the two sides of \eqref{eq:comparison} with the same signed
relative class.
\end{proof}

Therefore, one may use logarithmic invariants on the right-hand side
of \cite{WY}*{Theorem 1.3}. The comparison here includes the
$X_0$ insertion and the pushforward to $\scrC$, whereas an equality of
the unpushed virtual classes would be a different statement.
No argument for Theorem~\ref{thm:main} uses any type of orbifold
identity.

\begin{remark}
Identity \eqref{eq:comparison} is not a direct consequence
of \cite{BNR22}*{Theorem X}. Their theorem compares logarithmic
invariants of a pair with orbifold invariants of a
slope-sensitive modification of the pair, whereas \eqref{eq:comparison}
compares the two theories without modifying the target
$(P,X_\infty+X_\sigma)$, but with certain prescribed insertion and
pushforward. The additional input here is the degeneration
calculation for $(V,E)$, which, combined with
\cite{WY}*{Theorem 1.3}, establishes the desired comparison.
\end{remark}

\subsection{The comparisons used in the proof}\label{subsec:strategy}
We use the logarithmic theory with expansions of \cite{Ran}
whose virtual cycle pushes forward to the one for basic logarithmic maps given by
\cite{Ran}*{Theorem A and Proposition 4.6.1}\footnote{References to numbered sections of \cite{Ran} use the published
version, whose numbering differs from the arXiv version.}.
We distinguish unexpanded logarithmic moduli spaces, denoted by $\mathscr M$
from Ranganathan's moduli spaces of stable logarithmic maps with expansions, denoted by $\mathscr R$.
The latter are the moduli spaces used in the degeneration calculation. In our setting, the degeneration formulas given respectively by \cite{Ran} and \cite{KLR} will result in the same type of calculations regarding the vertex moduli spaces and gluing coefficients for a rigid tropical map, thus will give us the same conclusion at the end of the day. However, later in the paper, we will stick with Ranganathan's degeneration formula \cite{Ran} since the relevant perfect relative obstruction theory of the moduli stack of stable logarithmic maps with expansions is built over the corresponding Artin fan, which will give us a clean way for the comparison between the obstruction theories in question in order for us to claim that the corresponding virtual cycles match up.
For disjoint labeled boundary components, we permit independent
expansions.

For a smooth connected boundary, the smooth-pair comparison
\cite{AMW}*{Theorem 1.1.2} identifies these pushforwards with
Jun Li's relative theory. The distinction between expansion
stacks for a disconnected divisor is explained in
\cite{ACFW}*{Section 5.4.2} and no identification of those stacks
is used here. See Section~\ref{subsec:relative-degeneration}
for the precise degeneration formula and the vertex convention.

Write
\[
R=(\pi_{\log})_*[\mathscr M_X]^\vir,\qquad
L=(u_{\log})_*\left(\ev_0^*[X_0]\cap[\mathscr M_P]^\vir\right).
\]
Let $\mathscr R_X$ be the moduli space of stable logarithmic maps with expansions to $(X,D)$
with the same curve class, contact orders, and labels as $\mathscr M_X$, and set
\[
R_{\mathrm{rel}}=(\pi_{\mathrm{rel}})_*[\mathscr R_X]^\vir
\quad\text{in }A_*(\scrC)_{\QQ}
\] where $\pi_{\mathrm{rel}}$ denotes the natural projection to $\mathscr C$.

A smooth log--relative comparison gives $R=R_{\mathrm{rel}}$
by \cite{Ran}*{Proposition 4.6.1}, with a comparison to
Li's relative theory given by \cite{AMW}*{Theorem 1.1.2}.
In Section~\ref{sec:corner}, we construct the smooth pair $(V,E)$,
the lifted curve class $\lambda$, and the expanded logarithmic moduli space
\[
L_V=(u_V)_*\left(\ev_0^*[X_0]\cap[\mathscr R_V]^\vir\right).
\]
The proof is split into the following identities:
\begin{equation}\label{eq:proof-chain}
\underbrace{R=R_{\mathrm{rel}}}_{\text{smooth-pair comparison}},
\qquad
\underbrace{L=L_V}_{\text{Section~\ref{sec:corner}}},
\qquad
\underbrace{L_V=(-1)^{k_0}R_{\mathrm{rel}}}_{\text{degeneration formula calculation}}
\end{equation}
Only the last identity is the new relative calculation.
All these identities concern the prescribed data, insertion and
pushforward to $\scrC$. The paper is structured as follows.

\begin{enumerate}
\item Section~\ref{sec:corner} proves $L=L_V$. We blow up
$X_\infty\cap X_\sigma$ and a contact order $(k_i,k_i)$ becomes a
contact of order $k_i$ with the exceptional divisor. The pencil
spanned by $X_\infty$ and $X_\sigma$ shows that the other two
boundary divisors do not meet the maps in question. We then
explicitly switch from basic stable logarithmic maps to expanded logarithmic maps using the smooth-pair comparison theorem for $(V,E)$.

\item Section~\ref{sec:degeneration} constructs a degeneration
of the pair $(V,E)$. We first degenerate the base pair
$(X,D)$ and then perform the projective-bundle and corner-blowup
constructions over the whole family. This results in a degeneration of $(V,E)$ such that the general fiber is $(V,E)$ and
the central components are $Q=X\times\PP^1$ and $W$ over $D$.
The relative theory of $(X,D)$ will be obtained from the relative theory of $(X\times \PP^1, D\times \PP^1)$.

\item Sections~\ref{sec:vanishing} and~\ref{sec:caps} compute
which connected relative vertices on $W$ can contribute after
pushing forward to $\scrC$.
Lemma~\ref{lem:vanishing} eliminates those whose projections
give stable maps to $D$. The remaining contributions are
relative invariants of curves in a fiber, computed in
Section~\ref{sec:caps}.

\item In Section~\ref{sec:assembly}, we show that the only
remaining graph is a star shape. Each of $p_0,\ldots,p_n$ lies on its
own vertex over $D$. Including edge multiplicities, the
vertices carrying $p_i$ for $i>0$ contribute $1$, and the
vertex carrying $p_0$ contributes $(-1)^{k_0}$.
The central vertex gives $R_{\mathrm{rel}}$, proving the last
identity in \eqref{eq:proof-chain}. Combining the three identities
then proves the logarithmic theorem.
\end{enumerate}
Section~\ref{sec:vanishing} contains the full vanishing argument as to what type of graphs non-trivially contribute to the degeneration formula.
Constant projected maps are excluded by a standard dimension
count argument. For non-constant projections, a comparison of virtual
pullbacks reduces the assertion to the same dimension count
over a smooth ambient mapping space. Ordinary refined-diagonal
compatibility then gives the vanishing of the glued contribution.
We first assume $D$ to be connected. The cases $K=0$ and disconnected $D$
are addressed in Section~\ref{sec:endpoints}.

\section{From the SNC logarithmic target to a smooth pair}
\label{sec:corner}

The purpose of this section is to prove $L=L_V$.
The comparison with the relative/logarithmic theory of $(X,D)$ is deferred
to the degeneration calculation in the subsequent sections.

The contact orders on the projective bundle have the form $(k_i,k_i)$.
Blowing up the intersection of the two divisor components records
this common order as a contact order with the exceptional divisor.
This is the reason for considering
\[
b:V=\Bl_{X_\infty\cap X_\sigma}P\longrightarrow P.
\]
Let $E$ be its exceptional divisor. Write
$\widetilde X_\infty,\widetilde X_\sigma$ for the strict transforms of $X_\infty$ and $X_\sigma$ respectively.
The restriction of $\pi\circ b$ to $E$ gives $E\to D$ and let $e$ be its
fiber class and later we will abuse the notation $e$ for the pushforward class on $V$ along the inclusion $E\hookrightarrow V$. Let $f_V$ be the class of the strict transform of
a fiber of $P\to X$ over a point outside $D$, and put
$\ell=f_V-e \in H_2(V,\mathbb Z)$.
The section $X_0$ misses the blowup center, so we also abuse the notation $X_0$
 for its strict transform on $V$ and identify the curve class $\beta\in H_2(X)$
with its pushforward along this lifted section.

Notice that the pencil spanned by $X_\infty$ and $X_\sigma$ defines a rational map $P\dashrightarrow \mathbb P^1$, which is resolved to be a morphism by the blow-up operation, denoted by $g: V \longrightarrow \mathbb P^1$ where one has $g^*\OO_{\PP^1}(1)=b^*\OO_P(X_\infty)\otimes\OO_V(-E)$. In addition, it is an immediate fact that $g$ is constant in $X_0\subset V$, that is, $g(X_0)=c_0$ for some fixed $c_0\in\mathbb P^1$, and the fiber $g^{-1}(c_0)$ is disjoint from $\widetilde X_\infty\cup\widetilde X_\sigma$. 

\begin{lemma}\label{lem:pencil}
Every basic stable logarithmic map to
$(V,E+\widetilde X_\infty+\widetilde X_\sigma)$ whose composition
with $b$ has the curve class $\gamma$ and the labeled contact orders
prescribed for $\mathscr M_P$ has the curve class
\begin{equation}\label{eq:lambda}
\lambda=\beta+K\ell
\end{equation}
such that $b_*\lambda=\gamma$,
which satisfies
\[
E\cdot\lambda=K,\qquad
\widetilde X_\infty\cdot\lambda=\widetilde X_\sigma\cdot\lambda=0,\qquad
X_0\cdot\lambda=-k_0
\]

Moreover, 
every connected stable map of class $\lambda$ has an image fully landing in the fiber which contains $X_0$, and hence misses $\widetilde X_\infty\cup \widetilde X_\sigma$.
\end{lemma}

\begin{proof}
For a marking on a lifted logarithmic map of class $\lambda$ to
$(V,E+\widetilde X_\infty+\widetilde X_\sigma)$, let
$m_E,m_\infty,m_\sigma\in\mathbb Z_{\geq0}$ denote the
respective contact orders with
$E,\widetilde X_\infty,\widetilde X_\sigma$ at the marking.
The contact order with a boundary
component not containing the image of the marking is zero.
The divisor identities
\[
b^*X_\infty=\widetilde X_\infty+E,\qquad
b^*X_\sigma=\widetilde X_\sigma+E
\]
give the contact vector 
\[
(m_E+m_\infty,\;m_E+m_\sigma)
\]
after composition with $b$, ordered as
$(X_\infty,X_\sigma)$ by the functoriality of logarithmic maps.
Since $\widetilde X_\infty$ and $\widetilde X_\sigma$ are
disjoint, at least one of $m_\infty,m_\sigma$ is zero.
Therefore, for markings $p_i$, for $i>0$, the prescribed contact vector $(k_i,k_i)$
forces
\[
m_E=k_i,\qquad m_\infty=m_\sigma=0.
\]
provided that $k_i\neq 0$. So, we have $\widetilde X_\infty\cdot\lambda=\widetilde X_\sigma\cdot\lambda=0$. 
Similarly, a marking with contact vector $(0,0)$ downstairs
has all three contact orders zero upstairs. This suffices to determine the curve class $\lambda$ given by Equation (\ref{eq:lambda}).

Indeed, recall that $\lambda = \beta + K\ell$ and $b_*\lambda=\gamma$ where $\gamma=\beta+Kf$ and $\ell=f_V-e$. Here, we abuse the notation $\beta$ for the fixed curve class in $X$ and any possible lifting of it. Note that $b_*f_V=f$ since $b_*e=0$, thus $b_*(\lambda -\beta-Kf_V)=0$. The blowup formula implies that $\lambda=\beta+Kf_V+ae$ for some integer $a$ where only the coefficient $a$ needs to be determined. Furthermore, by projection formula, we know that 
\[
K = X_\infty\cdot\gamma = X_\infty\cdot b_*\lambda = b^*X_\infty\cdot \lambda = \widetilde X_\infty\cdot\lambda + E\cdot\lambda = E\cdot\lambda
\]
where $\widetilde X_\infty\cdot\lambda=0$ was obtained right above. 
Therefore, we have 
\[
\begin{aligned}
K = E\cdot\lambda
&=E\cdot(\beta+Kf_V+a e)\\
&=\underbrace{E\cdot\beta}_{0}
+K\underbrace{(E\cdot f_V)}_{0}
+a\underbrace{(E\cdot e)}_{-1}\\
&=-a,
\end{aligned}
\]
which concludes that $\lambda=\beta+kf_V-Ke=\beta+K(f_V-e)=\beta+K\ell$.
Finally, in $P$, recall that one has $X_0 = X_\infty - \pi^*D$ as a divisor class in $A^1(P)_\mathbb Q$. Therefore, we obtain 
\[
\begin{aligned}
X_0\cdot\lambda
&=b^*[X_0]\cdot\lambda\\
&=[X_0]\cdot b_*\lambda\\
&=([X_\infty]-\pi^*[D])\cdot\gamma\\
&=K-d\\
&=-k_0,
\end{aligned}
\]
which fully concludes the first half of the lemma.

For the second half, first of all, recall that the pencil spanned by divisors $X_\infty$ induces a rational map $P\dashrightarrow\mathbb P^1$ with base
locus $X_\infty\cap X_\sigma$. The blowup resolves this pencil
and gives a morphism $g: V\longrightarrow\mathbb P^1$. To describe its fibers,
choose a local equation $z=0$ for $D$ and the fiber
coordinates $[x:y]$ on $P$. The two sections $X_\infty$ and $X_\sigma$ have their respective defining 
equations $x=0$ and $x-zy=0$ while $X_0$ is given by $y=0$. So, locally, the morphism $g$ can be expressed via the formula $[x:x-zy]$, and its value on the entire $X_0$ is
a constant, denoted by $c_0=[1:1]$.

Let $[\xi_0:\xi_1]$ be the coordinates of the base $\mathbb P^1$.
Near $c_0$, set $w=1-\xi_1/\xi_0$. The closure of the graph
of the pencil has the local equation
\[
wx=zy.
\]
Indeed, the local expression of $g$ is given by $[x:y]\mapsto [\xi_0 : \xi_1]=[x:x-zy]$ and the local equation can be obtained by an easy algebra. 
Notice that this is the blowup of the ideal $(w,z)$ on $X\times\PP^1$, thus
away from $D\times\{c_0\}$, the map $(\pi\circ b,g)$ is an
isomorphism. Therefore, we actually have an isomorphism
\begin{equation}\label{eq:second-blowup}
V\simeq\Bl_{D\times\{c_0\}}(X\times\PP^1).
\end{equation}


Thus, there are two blow-down maps:
\[
\begin{tikzcd}
&V\arrow[dl,"b"']\arrow[dr,"{(\pi\circ b,g)}"]&\\
P&&X\times\PP^1.
\end{tikzcd}
\]
$X_0$ and $E$ are the strict transforms of
$X\times\{c_0\}$ and $D\times\PP^1$ respectively, along the second blowup map. Denoted by $F$ the exceptional divisor,
setting $w=0$ in the local equation gives $zy=0$, namely
the two components $F=\{w=0, z=0\}$ and $X_0=\{w=0, y=0\}$. Therefore, one has $g^{-1}(c_0) = X_0\cup F$. The divisors
$\widetilde X_\infty,\widetilde X_\sigma$ are two other fibers.

Recall that $g^*\OO_{\mathbb P^1}(1)=b^*\OO_P(X_\infty)\otimes\OO_V(-E)$, then nef line bundle $g^*\OO_{\mathbb P^1}(1)$ has degree
$X_\infty\cdot\gamma-E\cdot\lambda=K-K=0$ on $\lambda$.
Therefore, every component of the source has degree zero under
$g$. Since the source is connected, the whole image lies in
one fiber of $g$.
Also, 
\[
X_0\cdot\lambda=X_0\cdot(\beta+K\ell)=-d+K=-k_0<0,
\]
some non-constant component must lie in $X_0$. The pencil
value is therefore $c_0$, and the image misses
$\widetilde X_\infty\cup\widetilde X_\sigma$.
\end{proof}

Lemma~\ref{lem:pencil} only 
concerns basic stable logarithmic maps in the sense of Chen \cite{Chen},
Abramovich--Chen \cite{AC} and Gross--Siebert \cite{GS13}. For the degeneration formula and related calculations,
we use Ranganathan's moduli spaces of stable logarithmic maps with expansions \cite{Ran} following the
convention of Section~\ref{subsec:strategy}. The next proposition
makes this change of moduli spaces explicit through a comparison
of virtual classes with insertions.


\begin{proposition}
\label{prop:smooth}
Let $\mathscr R_V$ be the moduli stack of genus-zero stable logarithmic maps with expansions in the sense of \cite{Ran} to $(V,E)$ of class $\lambda$, with contact order $k_i$
at $p_i$ for $1\leq i\leq n$, and contact order zero with $p_0$ and the additional
$n_0$ interior markings. Its map
$u_V:\mathscr R_V\to\scrC$ records the projected stable map and
the $D$-evaluations of its $E$-markings. Then
\begin{equation}\label{eq:LV}
L=L_V:=(u_V)_*
\left(\ev_0^*[X_0]\cap[\mathscr R_V]^\vir\right).
\end{equation}
\end{proposition}
\begin{proof}
Logarithmic birational invariance in \cite{AW}*{Theorem 1.1.1}
applies to our construction here
\[
(V,\widetilde X_\infty+\widetilde X_\sigma+E)\longrightarrow(P,X_\infty+X_\sigma).
\]
Lemma~\ref{lem:pencil} gives the lifted data \eqref{eq:lambda}.
Moreover, by Lemma \ref{lem:pencil} again, these lifted maps do not meet
$\widetilde X_\infty$ or $\widetilde X_\sigma$.
We can therefore replace the target with the smooth pair $(V,E)$.
Let $\mathscr M_{V,E}^{\log}$ be the moduli space of basic stable
logarithmic maps to $(V,E)$ with these lifted data, and let
$u_{V,\log}:\mathscr M_{V,E}^{\log}\to\scrC$ be its projection.
The preceding argument gives
\[
L=(u_{V,\log})_*\left(\ev_0^*[X_0]\cap
[\mathscr M_{V,E}^{\log}]^\vir\right).
\]

We now pass to the logarithmic theory with expansions. Contracting the
expanded target, stabilizing the resulting map, and passing to
the basic logarithmic structure define the proper comparison morphism
\[
c_V:\mathscr R_V\longrightarrow\mathscr M_{V,E}^{\log}.
\]
By \cite{Ran}*{Proposition 4.6.1},
\[
(c_V)_*[\mathscr R_V]^\vir=[\mathscr M_{V,E}^{\log}]^\vir.
\]
Note that we have the commutative diagram 
\[
\begin{array}{ccc}
\mathscr R_V & \xrightarrow{c_V} & \mathscr M_{V,E}^{\log}\\
{\scriptstyle\operatorname{ev}_0\downarrow}
&&
\downarrow{\scriptstyle\operatorname{ev}_0}\\
V & \xrightarrow{\mathrm{id}} & V.
\end{array}
\]
Therefore, it follows from the projection formula that 
\[
\begin{aligned}
(c_V)_*\left(
(\operatorname{ev}_0)^*[X_0]
\cap[\mathscr R_V]^{\mathrm{vir}}
\right)
&=
(\operatorname{ev}_0)^*[X_0]
\cap(c_V)_*[\mathscr R_V]^{\mathrm{vir}}\\
&=
(\operatorname{ev}_0)^*[X_0]
\cap[\mathscr M_{V,E}^{\log}]^{\mathrm{vir}}.
\end{aligned}
\]
Hence, 
\[
\begin{aligned}
L_V
&=(u_V)_*\left(
\operatorname{ev}_0^*[X_0]\cap[\mathscr R_V]^{\mathrm{vir}}
\right)\\
&=(u_{V,\log})_*(c_V)_*\left(
\operatorname{ev}_0^*[X_0]\cap[\mathscr R_V]^{\mathrm{vir}}
\right)\\
&=(u_{V,\log})_*\left(
\operatorname{ev}_0^*[X_0]\cap
[\mathscr M_{V,E}^{\log}]^{\mathrm{vir}}
\right)\\
&=L.
\end{aligned}
\]

\end{proof}

\begin{remark}
The passage from basic stable logarithmic maps to Ranganathan's
logarithmic maps with expansions in Proposition~\ref{prop:smooth} prepares
for the degeneration argument in Section~\ref{sec:degeneration}.
There we apply the degeneration formula of
\cite{Ran}*{Theorem B and Section 7} to a degeneration of the
pair $(V,E)$ and express its virtual class through 
vertex moduli spaces and their gluing along the common boundary.

Proposition~\ref{prop:smooth} guarantees that such a switch of moduli
spaces preserves the enumerative problem
with the insertion $\ev_0^*[X_0]$.
We may therefore apply the degeneration formula \cite{Ran} using moduli spaces $\mathscr R$ of stable logarithmic maps with expansions,
and then interpret the comparison result in the moduli theory of
basic stable logarithmic maps in the sense of \cite{Chen}, \cite{AC} and \cite{GS13}.
\end{remark}

\section{Degenerating the resolved relative pair}\label{sec:degeneration}

From now on through Section~\ref{sec:assembly}, we aim to relate the pushforward virtual cycle $L_V$ to $R_{\operatorname{rel}}$ on $\mathscr C$ via the degeneration formula developed in \cite{Ran}. The core proposition to
be proved is as follows.
\begin{proposition}\label{prop:relative}
For the pushforward virtual cycle $L_V$ in Proposition~\ref{prop:smooth} and $R_{\mathrm{rev}}$ in Section \ref{subsec:strategy}, there is an identity of classes
\begin{equation}\label{eq:relative-calculation}
L_V=(-1)^{k_0}R_{\mathrm{rel}}\quad\text{in }A_*(\scrC)_{\QQ}.
\end{equation}
\end{proposition}
The proof will be completed in Section~\ref{sec:assembly}.
Together with the two comparison identities already established
in \eqref{eq:proof-chain}, this implies Theorem~\ref{thm:main}.

\subsection{Constructing the family of pairs $(V, E)$}
We construct a degeneration of $(V,E)$ by first degenerating
its base pair $(X,D)$ and then repeating the bundle and corner-blowup
construction over the whole family. Let $t$ denote the degeneration
parameter and it is distinct from the coordinate on the pencil base 
$\PP^1$ used in Section~\ref{sec:corner}.
Consider the deformation to the normal cone
\[
\mathfrak X=\Bl_{D\times\{0\}}(X\times\AA^1),\qquad
\mathfrak D=\widetilde{D\times\AA^1}.
\]
Its generic fiber is $(X, D)$ and the central fiber is $X\cup_D Y$, where, writing
$N_D=\OO_D(D)$,
\[
Y=\PP_D(N_D^{-1}\oplus\OO_D).
\]
Denote the zero section of $Y$ by $D_0$ and the infinity
section by $D_\infty$. The gluing takes place along $D_0$ and $\mathfrak D$ specializes to $D_\infty$, that is, $\mathfrak D\cap\mathfrak X_0=D_\infty$.

Now set
\[
\mathfrak P=\PP_{\mathfrak X}
(\OO_{\mathfrak X}(-\mathfrak D)\oplus\OO).
\]

Since $\mathfrak D\cap X=\varnothing$ inside $\mathfrak X_0$, the
projective bundle is just a trivial bundle $X\times\PP^1$ over $X$.
Analogously, the section $\sigma_{\mathfrak X}$ defining $\mathfrak D$ defines the 
section $\mathfrak X_\sigma$ of $\mathfrak P$, just as $\sigma$ defines those $X_\sigma$ of $P=\mathbb P(\mathcal O_X(-D)\oplus \mathcal O_X)$.
Blow up the intersection of $\mathfrak X_\sigma$ and $\mathfrak X_\infty$ and write
$\mathfrak V\to\AA^1$ for the resulting family and
$E_{\mathfrak V}$ for its exceptional divisor. To sum up,
for every $t\ne0$, the generic fibers of this chain of constructions are respectively
\[
(\mathfrak X_t,\mathfrak D_t)=(X,D),\qquad
\mathfrak P_t=P,\qquad
(\mathfrak V_t,E_{\mathfrak V,t})=(V,E).
\]
Hence, the family $\mathfrak V\to\AA^1$ degenerates exactly to the relative target of
Proposition~\ref{prop:smooth}.
To see the central fiber $\mathfrak V_0$ of $\mathfrak V\rightarrow\mathbb A^1$ and that no extra central component is created during this entire process, we use the quotient
coordinates locally (following our projectivization convention). Let $z$ be a local defining equation of
$\mathfrak D$, and let $[x:y]$ be the quotient fiber coordinates
given by the two summands in compatible local trivializations.
The two sections $\mathfrak X_\infty$ and $\mathfrak X_\sigma$ have their respective local equations
\[
x=0,\qquad x-zy=0,
\]
and the zero section is locally given by $y=0$. The intersection of $\mathfrak X_\infty$ and $\mathfrak X_\sigma$ has the ideal
$(x,z)$ in the chart where $y$ is invertible. It is smooth
over $\AA^1$, lies over $\mathfrak D$, and misses the gluing divisor $D$
of $\mathfrak X_0$. 

Therefore, 
the central fiber $\mathfrak V_0$ consists of $Q=X\times\PP^1$, and the variety
$W$ obtained by applying the projective-bundle and
corner-blowup construction to $(Y,D_\infty)$, that is, $W=\operatorname{Bl}_{s_\infty(D_\infty)}\mathbb P_Y\!\left(\mathcal O_Y(-D_\infty)\oplus\mathcal O_Y\right)$ where $s_\infty:Y\rightarrow\mathbb P_Y\!\left(\mathcal O_Y(-D_\infty)\oplus\mathcal O_Y\right)$ is the section induced by the quotient onto $\mathcal O_Y$. It is obvious that 
$Q$ and $W$ meet in $F_0=D\times\PP^1$.
We will use $F_0$ for both copies of this gluing divisor and the inverse image of $D_0$ on $W$.
Note that the central-fiber restriction of $E_{\mathfrak V}$ is supported on $W$.
Write $E_W=E_{\mathfrak V}|_W$, then $E_W\cap F_0=\varnothing$.
The resulting central fiber pair is
\begin{equation}\label{eq:central}
(Q,F_0)\ \cup_{F_0}\ (W,F_0+E_W).
\end{equation}
Here $E_W$ is the specialization of the original relative boundary
$E$, whereas $F_0$ is the new gluing boundary on both central
components. In particular, the retained markings $p_i$, $i>0$,
are relative to $E_W$ and lie on the $W$-side. The gluing markings
relative to $F_0$ are additional labels supplied by the degeneration
formula. The interior marking $p_0$ may initially lie on either side.

Furthermore, let $g_{\mathfrak V}:\mathfrak V\to\PP^1$ be the resolved
family pencil analogous to the one in Section \ref{sec:corner}, given locally by $[x:x-zy]$.
It restricts to $g:V\rightarrow \mathbb P^1$ on the general fiber and to an automorphism
of the $\PP^1$ factor on $Q$.
Let $X_{0,\mathfrak V}$ be the strict transform of the zero section of the total projective bundle $\mathfrak P$. It extends $X_0$ and likewise 
has pencil value $c_0=[1:1]$.
Write $X_{0,Q}=X_{0,\mathfrak V}|_Q=X\times\{c_0\}$ and
$X_{0,W}=X_{0,\mathfrak V}|_W$.
The restriction $g_W=g_{\mathfrak V}|_W$ is the pencil on $W$.
Since the total degree along the pencil is zero
and the pencil line is nef, every effective class assigned to a vertex in a
degeneration splitting has degree zero along the pencil.

\subsection{The curve classes in the degeneration}\label{subsec:classes}
Since we will apply the degeneration formula with the fixed curve class $\lambda$ in a general fiber, the following discussion is dedicated to laying down several crucial divisor class identities on the central fiber during the degeneration we constructed, which will be used in the degeneration calculations later. 

Let $\pr_X:\mathfrak V\to X$ be the composition of the bundle
projection, the blow-down to $X\times\AA^1$, and the projection to $X$.
Let $h\in A^1(\mathfrak V)$ be the class of $\widetilde{\mathfrak X}_\infty$ (same as that of $\widetilde{\mathfrak X}_\sigma$), which is the pullback of $\mathfrak X_\infty$ to the blowup $\mathfrak V\rightarrow\mathfrak P$.
Let \(Y_{\mathrm{exc}}\) denote the pullback to \(\mathfrak V\) of the exceptional divisor \(Y\subset\mathfrak X\) of the blow-up \(\mathfrak X= \Bl_{D\times \{0\}}(X\times\mathbb A^1)\to X\times\mathbb A^1\).
The following divisor class identities hold on the total space $\mathfrak V$ with $h=c_1(\mathcal O_{\mathbb P^1}(1))\in A^1(\mathbb P^1)$:
\begin{equation}\label{eq:global-relations}
\begin{split}
[X_{0,\mathfrak V}]&=h-\pr_X^*[D]+[Y_{\mathrm{exc}}],\\
g_{\mathfrak V}^*H&=h-[E_{\mathfrak V}].
\end{split}
\end{equation}
Indeed, let $q:\mathfrak V\rightarrow\mathfrak X$ be the composition of the blow-down map $\mathfrak V\rightarrow\mathfrak P$ and projective bundle projection $\mathfrak P\rightarrow\mathfrak X$. Notice that $\text{total transform of }D\times\mathbb A^1
=\mathfrak D+Y$, and pulling the equality back to $\mathfrak V$ gives rise to $\operatorname{pr}_X^*[D]
=q^*[\mathfrak D]+[Y_{\mathrm{exc}}]$. Moreover, the projective bundle yields the relation $[X_{0,\mathfrak V}]=h-q^*[\mathfrak D]$, and thus one obtains 
\[
[X_{0,\mathfrak V}]
=h-\bigl(\operatorname{pr}_X^*[D]-[Y_{\mathrm{exc}}]\bigr)
=h-\operatorname{pr}_X^*[D]+[Y_{\mathrm{exc}}].
\]
The second equality of (\ref{eq:global-relations}) is just due to the construction of the pencil $g_\mathfrak V$.

The restrictions of $[Y_{\mathrm{exc}}]$ to $Q$ and $W$ are $[F_0]$
and $-[F_0]$, respectively. Each edge has the same contact order
on its two sides. Its two contributions to the degree of
$Y_{\mathrm{exc}}$ therefore cancel.
Indeed, the two blowups and the projective-bundle formula give,
for connected $D$,
\[
H_2(\mathfrak V,\mathbb Z)\simeq H_2(X,\mathbb Z)\oplus\mathbb Z^3.
\]
The three additional bases are 
$[Y_{\mathrm{exc}}],h$, and $[E_{\mathfrak V}]$ and the general-fiber
class $\lambda$ has respective intersection numbers $0,K$, and $K$ with them.
Contact orders along gluing edges have total \(Y_{\mathrm{exc}}\)-degree zero in the central fiber, and the total intersections against
the other two divisors are prescribed to be $K$.
These conditions are sufficient to determine the total class of the
specialization of $\lambda$. Nonetheless, applying the degeneration formula still
requires summing over all compatible choices of vertex classes with a fixed total class. 

\subsection{The degeneration formula with expansions}
\label{subsec:relative-degeneration}

Let $\Gamma$ be a connected bipartite tree with vertices over
$Q$ and $W$, decorated by effective curve classes,
relative and interior markings, and contact orders at the relative markings as in Section~\ref{subsec:classes}.
For each vertex $v$, let $\mathscr R_v$ be the moduli space of genue-zero stable logarithmic maps with expansions to $(Q,F_0)$ or $(W,F_0+E_W)$
with the data prescribed at $v$ using the convention of
Section~\ref{subsec:strategy}.
Denote by $a_e$ the contact order at an edge $e$.

Let $\Delta_\Gamma$ be the product of the diagonals of $F_0$,
one for each edge. Define $\mathcal G_\Gamma$ by the Cartesian
square
\[
\begin{tikzcd}[column sep=large]
\scrC
&
\mathcal G_\Gamma
  \arrow[l,"u_\Gamma"']
  \arrow[r] \arrow[d]
&
\displaystyle\prod_v\mathscr R_v
  \arrow[d,"\ev"]
\\
&
F_0^{E(\Gamma)}
  \arrow[r,"\Delta_\Gamma"']
&
(F_0\times F_0)^{E(\Gamma)}.
\end{tikzcd}
\]
Here, $\ev$ records the evaluations at the gluing markings,
and $u_\Gamma$ is the natural forgetful map.
We also write
$\ev_0:\prod_v\cR_v\to\mathfrak V$
for evaluation at $p_0$, followed by contraction to the
original target and inclusion into the total space.

\begin{proposition}\label{prop:relative-degeneration}
There is a decomposition formula
\begin{equation}\label{eq:relative-degeneration-explicit}
L_V=
\sum_\Gamma
\frac{\prod_{e\in E(\Gamma)}a_e}
     {|\operatorname{Aut}(\Gamma)|}
(u_\Gamma)_*\Delta_\Gamma^!
\left(
\ev_0^*[X_{0,\mathfrak V}]
\cap\boxtimes_{v\in V(\Gamma)}[\cR_v]^\vir
\right)
\end{equation}
in $A_*(\scrC)_\QQ$.
Here, $\boxtimes$ denotes the exterior product of cycles.
An automorphism preserves the graph, curve classes, and marking profiles with contact orders.
\end{proposition}

\begin{proof}
We will apply \cite{Ran}*{Theorem B} to the family $(\mathfrak V, \mathfrak V_0 + E_{\mathfrak V})\rightarrow(\mathbb A^1, 0)$. First of all, by our construction of this family, which is discussed in Section \ref{sec:degeneration}, it is obvious that the family is toroidal degeneration without self-intersection. 

Let $\mathscr K_\rho:=K_\rho(\mathfrak V_0)$ be the Ranganathan's moduli space of stable logarithmic maps with expansions marked by a rigid splitting type $\rho$, which is a bipartite tree. By \cite{Ran}*{Theorem B}, we have the following decomposition 
\[
[\mathscr R_0]^{\mathrm{vir}}
=
\sum_\rho m_\rho\,(\iota_\rho)_*
[\mathscr K_\rho]^{\mathrm{vir}}.
\]
Note that we have the following diagrams
\[
\begin{array}{ccccc}
\mathscr K_\rho
&\xrightarrow{\ \mu_\rho\ }&
\mathcal G_\rho
&\longrightarrow&
\mathscr P_\rho
\\[6pt]
{\scriptstyle u_0\circ\iota_\rho}\downarrow
&&
\downarrow{\scriptstyle\phi_\rho'}
&&
\downarrow{\scriptstyle\phi_\rho}
\\[6pt]
\mathscr C
&\xleftarrow{\ u_\Gamma\ }&
\mathcal G_\Gamma
&\longrightarrow&
\displaystyle\prod_v\mathscr R_v
\\[6pt]
&&\downarrow&&\downarrow{\scriptstyle\mathrm{ev}}
\\[6pt]
&&F_0^{E(\Gamma)}
&\xrightarrow{\ \Delta_\Gamma\ }&
(F_0\times F_0)^{E(\Gamma)}.
\end{array}
\]
where $
\mathscr P_\rho
:=
\mathop{\times}\limits_v^\lambda K_\rho(X_v)
$ is introduced in Part C of \cite{Ran}*{Theorem B} and $\mathcal G_\rho
=
\mathscr P_\rho
\times_{(F_0\times F_0)^{E(\Gamma)}}
F_0^{E(\Gamma)}$.
Ranganathan's gluing formula yields the virtual class identities
\[
[\mathcal G_\rho]^{\mathrm{vir}}
=\Delta_\Gamma^![\mathscr P_\rho]^{\mathrm{vir}},
\qquad
[\mathscr K_\rho]^{\mathrm{vir}}
=\mu_\rho^*[\mathcal G_\rho]^{\mathrm{vir}}.
\]
Hence, one obtains 
\[
[\mathscr K_\rho]^{\mathrm{vir}}
=
\mu_\rho^*\Delta_\Gamma^!
[\mathscr P_\rho]^{\mathrm{vir}}.
\]
Furthermore, since $\mu_\rho$ is a finite \'etale map and let \(d_\rho=\deg\mu_\rho\), the Cartesian diagram gives rise to 
\[
\begin{aligned}
(\phi_\rho'\circ\mu_\rho)_*
[\mathscr K_\rho]^{\mathrm{vir}}
&=
d_\rho(\phi_\rho')_*
\Delta_\Gamma^![\mathscr P_\rho]^{\mathrm{vir}}\\
&=
d_\rho\Delta_\Gamma^!
(\phi_\rho)_*[\mathscr P_\rho]^{\mathrm{vir}}\\
&=
\boxed{
d_\rho\Delta_\Gamma^!
\left(\boxtimes_v[\mathscr R_v]^{\mathrm{vir}}\right)
}
\end{aligned}
\]
So, we can express the moduli space of stable logarithmic maps with expansions to the central fiber $\mathscr R_0$ as 
\[
[\mathscr R_0]^{\mathrm{vir}}
=
\sum_\rho m_\rho(\iota_\rho)_*
[\mathscr K_\rho]^{\mathrm{vir}}
=
\sum_\rho m_\rho(\iota_\rho)_*\mu_\rho^*
\Delta_\Gamma^![\mathscr P_\rho]^{\mathrm{vir}}.
\]
Notice that one has $
u_0\circ\iota_\rho
=
u_\Gamma\circ\phi_\rho'\circ\mu_\rho
$, and then pushing the decomposition above to $\mathscr C$ gives the identity
\[
(u_0)_*[\mathscr R_0]^{\mathrm{vir}}
=
\sum_\rho m_\rho d_\rho\,
(u_\Gamma)_*\Delta_\Gamma^!
\left(\boxtimes_v[\mathscr R_v]^{\mathrm{vir}}\right)
=
\sum_\Gamma\frac{\prod_ea_e}{|\operatorname{Aut}(\Gamma)|}(u_\Gamma)_*\Delta_\Gamma^!
\left(\boxtimes_v[\mathscr R_v]^{\mathrm{vir}}\right)
\]
where $m_\rho d_\rho=\frac{\prod_ea_e}{|\operatorname{Aut}(\Gamma)|}$ and this factor is also discussed in \cite{KLR}.

Together with the projection formula and the compatibility of cap products with refined Gysin pullback, we have the desired decomposition formula. 
\end{proof}

Next, we need to identify rigid splitting graphs which stratify the total moduli stack in \cite{Ran}*{Theorem B} for the degeneration. We claim that the underlying graph of every rigid splitting type is a bipartite tree with vertices labeled with $Q$ and $W$.

\begin{lemma}\label{bipartite}
    Every connected genus-zero rigid tropical splitting type for the degeneration $\mathfrak V\longrightarrow\mathbb A^1$ has a bipartite underlying graph.
\end{lemma}
\begin{proof}
Indeed, tropicalizing the family results in the following tropical map
\[
\Sigma_{\mathfrak B}\longrightarrow\mathbb R_{\geq0},
\qquad
(a,b,c)\mapsto a+b.
\]
where the tropicalization $\Sigma_{\mathfrak V}$ of $\mathfrak V$ can be expressed as 
\[
\Sigma_{\mathfrak B}
=
\{(a,b,0):a,b\geq0\}
\ \cup\
\{(0,b,c):b,c\geq0\}.
\]
and the tropicalization of the base $(\mathbb A^1, 0)$ is just $\mathbb R_{\geq 0}$.
Let us take the slice of $\Sigma_{\mathfrak V}$ over $1\in\mathbb R_{\geq 0}$, denoted by $\Sigma_1$, that is, $
\Sigma_1
=
\{(a,b,0):a+b=1\}
\cup
\{(0,1,c):c\geq0\}.
$
Let $Q=(1, 0, 0), W=(0, 1, 0)\in\Sigma_1$. Here, we abuse the notation $Q$ and $W$ for the respective tropical point in $\Sigma_1$ corresponding to the actual components $Q$ and $W$ of $\mathfrak V_0$.
For any tropical map $h:\Gamma\longrightarrow \Sigma_1$, if there is a vertex $v$ such that $h(v)\in[Q, W]^\circ$, then such a tropical map has one free length parameter making it admit a continuous deformation, which is not rigid. Hence, a rigid splitting type in our scenario must be a bipartite tree. 

\begin{figure}[htbp]
  \centering
  \includegraphics[width=0.8\textwidth]{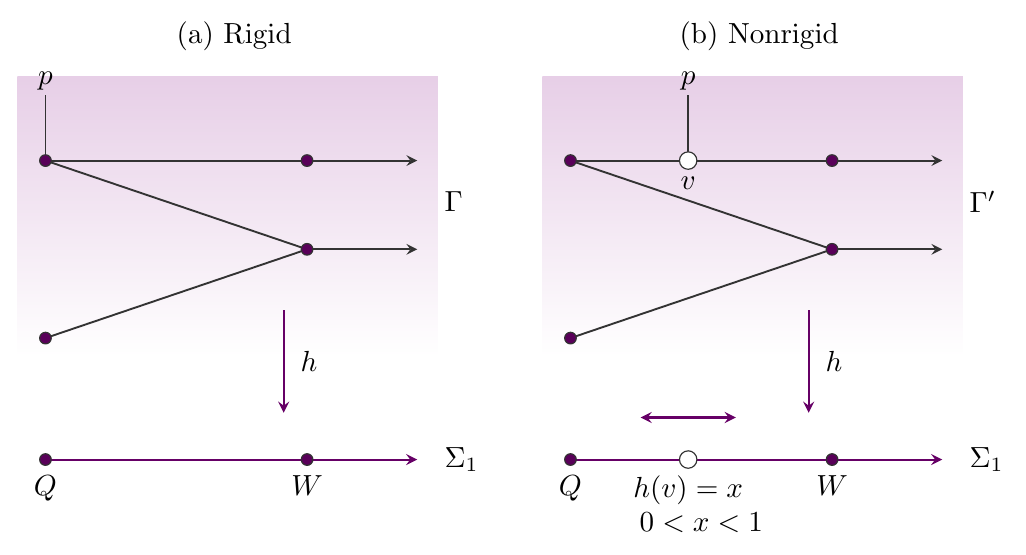}
  \caption{If $h(v)=x$ sits in the interior of the line segment $[Q, W]$, then $x$ is a free parameter resulting a one-dimensional continuous deformation of $h$, which is not a rigid tropical map.}
  \label{fig:rigid-splitting}
\end{figure}
\end{proof}

\section{Vanishing for vertices over \texorpdfstring{$W$}{W}}
\label{sec:vanishing}

We now analyze the contributions of the vertex moduli spaces $\mathscr R_v$ in the degeneration formula of $L_V$ for a rigid tropical type when $v$ belongs to $W$-vertices. Throughout this analysis, by a vertex moduli space having ``vanishing
contribution'', we mean that the corresponding virtual cycle vanishes after gluing with another factors and pushing forward to $\scrC$. It does not necessarily indicate 
that the virtual cycle of the moduli space itself is zero.

Prior to any vanishing result, we need to settle down some notation. Let $\rho:W\to D$ be the projection and put $h_W=h|_W$.
Fix the data of a connected $W$-vertex including the curve class
$\delta\in H_2(W,\mathbb Z)$, the number of markings denoted by $q$, and the gluing flags.
Also,  write $\delta_D=\rho_*\delta$ for the projected curve class.
The ordered contact profiles at $F_0$ and $E_W$ are part of
these fixed data. A pencil-point insertion at a flag means
$\ev^*g_W^*H$ at that flag.

The purpose of this
section is to show that for a moduli space associated with a $W$-vertex with fixed curve class $\delta$ and number $q$ of markings, it has a non-trivial contribution to the decomposition of $L_V$ only if $\delta_D=0$ and $q\leq 2$, meaning that
the stable logarithmic maps parametrized by the moduli space become constant after projection
to $D$ and it has at most two markings, including
the gluing markings. In that case, we will only need to study the logarithmic Gromov-Witten theory 
of a fiber surface pair $
\bigl(S_d:= \rho^{-1}(d),\,(F_0\cap S_d)+(E_W\cap S_d)\bigr)
$ for any point $d\in D$.

\begin{lemma}[Logarithmic first Chern class]\label{lem:c1}
Given the notation as above, the following identity holds
\begin{equation}\label{eq:c1}
c_1\bigl(T_W(-\log(F_0+E_W))\bigr)
=\rho^*c_1(T_D)+2g_W^*H.
\end{equation}
Consequently, when the projected data $(\delta_D,q)$ are stable
and $g_{W,*}\delta=0$, one obtains 
\[
\operatorname{vdim}\mathscr R_v = \operatorname{vdim}\mathscr M_{0,q}(D,\delta_D) + 2.
\]
\end{lemma}
\begin{proof}
First of all, recall that we have the following chain of maps in the construction of $W$ with $\rho:=p\circ r\circ b_W$
\[
W\xrightarrow{\ b_W\ }Z\xrightarrow{\ r\ }Y\xrightarrow{\ p\ }D
\]
where
\[
Y=\PP_D(N_D^{-1}\oplus\OO_D),\qquad
Z=\PP_Y(\OO_Y(-D_\infty)\oplus\OO_Y).
\]
Here $p$ and $r$ are the projections, and $b_W$
blows up the intersection of the two sections $s_\infty$ and $s_\sigma$ of $Z$ with the exceptional divisor $E_W$, where $s_\infty$ and $s_\sigma$ are analogously defined as before.
Let
$\eta=c_1(\OO_Y(1))$ and $\xi=c_1(\OO_Z(1))$.
With our quotient convention, the section classes and their pullbacks satisfy
\[
D_0=\eta,\qquad D_\infty=\eta+p^*c_1(N_D),\qquad
h_Z=\xi+r^*D_\infty,\qquad h_W=b_W^*h_Z,
\]
where $h_Z$ is the common class of the two positive sections of $Z$.
The gluing divisor is $F_0=b_W^*r^*D_0$ since the blowup
center lies over $D_\infty$ and is disjoint from $D_0$.
The resolved pencil $g_W:W\to\PP^1$ satisfies
$g_W^*H=h_W-E_W$, where $H=c_1(\OO_{\PP^1}(1))$.

Now we are in the position to calculate $c_1\bigl(T_W(-\log(F_0+E_W))\bigr)$ with the divisor class identities above. Recall that for a rank-two bundle $\mathcal E$ and the quotient projection
$\pi:\PP_B(\mathcal E)\to B$, the relative Euler sequence gives
\[
c_1(T_{\PP_B(\mathcal E)/B})
=2c_1(\OO_{\PP_B(\mathcal E)}(1))-\pi^*c_1(\mathcal E).
\]

First, via the two exact sequences of tangent bundles
\[
\begin{aligned}
0&\longrightarrow T_{Y/D}\longrightarrow T_Y
  \longrightarrow p^*T_D\longrightarrow 0,\\
0&\longrightarrow T_{Z/Y}\longrightarrow T_Z
  \longrightarrow r^*T_Y\longrightarrow 0,
\end{aligned}
\]
we deduce that \[
\begin{aligned}
    c_1(T_Z)
&=(p\circ r)^*c_1(T_D)
+r^*c_1(T_{Y/D})
+c_1(T_{Z/Y})\\
&=(p\circ r)^*c_1(T_D) + r^*(2\eta + p^*c_1(N_D)) + (2\xi+r^*D_\infty)\\
&=(p\circ r)^*c_1(T_D) + r^*(D_0 + D_\infty) + (2h_Z-r^*D_\infty)\\
&=(p\circ r)^*c_1(T_D)+r^*D_0 + 2h_Z
\end{aligned}
\]

Note that $
c_1(T_W)=b_W^*c_1(T_Z)-E_W.
$
Therefore, we have 
\[
\begin{aligned}
c_1(T_W) &= b_W^*(p\circ r)^*c_1(T_D)+b_W^*r^*D_0 + 2b_W^*h_Z-E_W\\
&=\rho^*c_1(T_D)+F_0+2h_W-E_W
\end{aligned}
\]
where the second equality is due to the fact that $b_W^*r^*D_0=F_0$ and $h_W=b_W^*h_Z$.

Hence, one obtains 
\[
\begin{aligned}
 c_1\!\left(T_W(-\log(F_0+E_W))\right) &= c_1(T_W)-F_0-E_W \\
 &= \rho^*c_1(T_D)+F_0+2h_W-E_W -F_0 - E_W\\
 &=\rho^*c_1(T_D) + 2(h_W-E_W)\\
 &= \rho^*c_1(T_D) + 2g_W^*H.
\end{aligned}
\]
where the last equality holds because $g_W^*H=h_W-E_W$. This is exactly the equality we demand. 

Lastly, we claim that $\text{vdim}\mathscr R_v=\text{vdim}\mathscr M_{0, q}(D, \delta_D) + 2$ when logarithmic maps to $D$ prescribed with ($\delta_D, q$) are stable and $g_{W,*}\delta = 0$. Indeed, at first, we notice that 
\[
\begin{aligned}
c_1\!\left(T_W(-\log(F_0+E_W))\right)\cdot\delta
&=c_1(T_D)\cdot\rho_*\delta+2H\cdot g_{W,*}\delta\\
&=c_1(T_D)\cdot\delta_D.
\end{aligned}
\]
Therefore, we have 
\[
\begin{aligned}
\operatorname{vdim}\mathscr R_v
&=\dim W-3+q+c_1\!\left(T_W(-\log(F_0+E_W))\right)\cdot\delta\\
&=\bigl(\dim D-3+q+c_1(T_D)\cdot\delta_D\bigr)+2\\
&=\operatorname{vdim}\mathscr M_{0,q}(D,\delta_D)+2.
\end{aligned}
\]
\end{proof}

The next lemma shows that the vertex moduli space with stable projected
data to $D$ contributes zero after gluing and pushforward to $\scrC$. For a vertex moduli space $\mathscr R_D$,
having two extra virtual dimensions on $\mathscr M_{0, q}(D, \delta_D)$ shown in Lemma~\ref{lem:c1}
suggests this vanishing. However, the moduli space of the maps
to $D$ need not have its expected dimension, so a
virtual-dimension count argument alone does not suffice to give the vanishing statement.

Prior to the statement of Lemma \ref{lem:vanishing}, one needs to settle down some notation. Recall that we have a chain of maps \[
\mathscr R_D\xrightarrow{\operatorname{ev}_j}F_0
\xrightarrow{\,g_W|_{F_0}\,}\mathbb P^1.
\]
Let $\zeta$ be defined as follows
\[
\zeta=
\begin{cases}
1, & \text{if there is no insertion},\\[2pt]
(g_W\circ\operatorname{ev}_j)^*H,
& \text{for a gluing marking }j\text{ along }F_0,\\[2pt]
\operatorname{ev}_0^*[X_{0,W}],
& \text{if the vertex carries }p_0.
\end{cases}
\]
The following lemma asserts that the virtual cycle $[\mathscr R_D]^{\operatorname{vir}}$ capping with $\zeta$ whichever case happens to $\zeta$ will push forward to $0$ as long as the projected data to $D$ are stable and $g_{W,*}\delta=0$.

\begin{lemma}\label{lem:vanishing}
Assume that $(\delta_D=\rho_*\delta,q)$ is stable, that is,
$\delta_D\ne0$ or $q\geq3$. Let $\mathscr R_D$ be the connected
vertex moduli stack with expansions of $(W,F_0+E_W)$ with
the fixed class, labels, and contact order profiles just specified.
Let
\[
p_D:\mathscr R_D\longrightarrow\mathscr B_D=\mathscr M_{0,q}(D,\delta_D)
\]
be the projection and stabilization. Suppose $g_{W,*}\delta=0$.
One has
\begin{equation}\label{eq:vanish}
(p_D)_*\left(\zeta\cap[\mathscr R_D]^\vir\right)=0.
\end{equation}
Consequently, every term of the degeneration formula containing
this inserted vertex vanishes after gluing and pushing forward to
$\scrC$.
\end{lemma}

\begin{proof}
Let $a=\operatorname{codim}\zeta\leq1$. First of all, if $\delta_D=0$, then
$q\geq3$ and hence $\mathscr B_D=D\times\overline{\mathscr M}_{0,q}$ is smooth of the
expected dimension, where $\overline{\mathscr M}_{0,q}$ is the moduli stack of
$q$-marked stable rational curves. Lemma~\ref{lem:c1} gives that
$\operatorname{vdim}\mathscr R_D-a>\dim\mathscr B_D$ which concludes
\eqref{eq:vanish} directly by a standard virtual dimension counting argument.

For $\delta_D\ne0$, $\mathscr B_D$ can be
obstructed and a straightforward virtual dimension count argument does not go through. The strategy here to get around a possible obstruction on $\mathscr B_D$ is to reduce to a smooth and convex ambient space, and
pull all the vanishing results back based on compatibility of obstruction theories.

Pick a very ample line bundle $M$ on $D$ such that $M\otimes N_D$ is very ample as well in which $N_D=\OO_D(D)$.
The corresponding complete linear systems define an embedding
$i:D\hookrightarrow T=\PP^{r_1}\times\PP^{r_2}$ with
\[
i^*\OO_T(1,0)\simeq M,\qquad
i^*\OO_T(0,1)\simeq M\otimes N_D.
\]
Hence $i^*\OO_T(-1,1)\simeq N_D$, and the line bundle
used to construct $W$ extends to $T$.
Construct $(W_T,F_{0,T}+E_T)$ over $T$ using $\OO_T(-1,1)$ in the same way as $W$ is constructed.
Therefore, one has the following Cartesian diagram
\[
\begin{array}{ccc}
W & \xrightarrow{\ j\ } & W_T\\
{\scriptstyle\rho}\downarrow && \downarrow{\scriptstyle\rho_T}\\
D & \xrightarrow{\ i\ } & T .
\end{array}
\]
So that the boundary divisors also restrict, that is, $
j^*F_{0,T}=F_0$ and $j^*E_T=E_W$. In particular, $\rho_T$ is smooth and $j$ is regular. Moreover, the resolved pencil $g_W:W\rightarrow\mathbb P^1$ also extends to $W_T$, meaning that there exists a pencil $g_T:W_T\rightarrow\mathbb P^1$ satisfying $g_T\circ j=g_W$ and $j^*[X_{0,W_T}]=[X_{0,W}]$. So, by the definition of $\zeta$, it extends to some $\zeta_T$ on $W_T$ of codimension $a$.

Let $\mathscr R_T$ be the corresponding moduli stack of logarithmic stable maps with expansions to $W_T$ having the class $j_*\delta$ and the same labeled
markings, and put $\mathscr B_T=\mathscr M_{0,q}(T,i_*\delta_D)$.
Projection and stabilization give the following Cartesian diagram
\begin{equation}\label{eq:vanishing-square}
\begin{tikzcd}
\mathscr R_D\arrow[r,"\jmath"]\arrow[d,"p_D"']&
\mathscr R_T\arrow[d,"p_T"]\\
\mathscr B_D\arrow[r,"\iota"']&\mathscr B_T.
\end{tikzcd}
\end{equation}

Now, we need to verify that $\iota$ defines a virtual pullback needed for this square. Since \cite{Manolache}*{Proposition 5.3} gives the relative obstruction theory
\[
\mathbb E_\iota
=\bigl(R\pi_{D*}\varphi_D^*N_{D/T}\bigr)^\vee[1]
\longrightarrow
\mathbb L_{\mathscr B_D/\mathscr B_T}
\]
for the map $\iota$, it suffices to verify that this obstruction is perfect, i.e. \(H^1(C,\varphi^*N_{D/T})=0\) at every point $[\varphi:C\rightarrow D]\in\mathscr B_D$. Indeed, since $N_{D/T}$ is the normal bundle of the chosen embedding $D\hookrightarrow T$ and $T$ is the product of projective spaces, it yields the tangent bundle $T_T$ of $T$ is globally generated, and so is $N_{D/T}$ due to the surjection \(i^*T_T\twoheadrightarrow N_{D/T}\). This global generation property gives rise to a short exact sequence with some kernel $K$
\[
0\longrightarrow K\longrightarrow\mathcal O_C^{\oplus m}
\longrightarrow\varphi^*N_{D/T}\longrightarrow0.
\]
Applying the long exact sequences in cohomology corresponding to the short exact sequence, one has \(H^1(C,\varphi^*N_{D/T})=0\). 
Writing $\pi_D,\varphi_D$ for the respective universal curve and map
over $\mathscr B_D$, the sheaf $\cV=\pi_{D*}\varphi_D^*N_{D/T}$ is represented by a single
vector bundle and commutes with base change. Then
\cite{Manolache}*{Proposition 5.3} gives the normal obstruction
theory $\mathbb E_\iota=\cV^\vee[1]$.

The target pairs $(W, F_0+E_W)$ and $(W_T, F_{0,T}+E_T)$ both have the same Artin fan
\[
\mathscr A=
[((\AA^1)^2\setminus\{(0,0)\})/(\mathbb G_m)^2]
\]
since $F_0\cap E_W=\varnothing$ and $F_{0,T}\cap E_T=\varnothing$.
Their maps to $\mathscr A$ are smooth, and $j$ is a morphism
over $\mathscr A$.

Over a finite-type open substack of the moduli stack of prestable logarithmic maps to $\mathscr A$, 
after choosing a subdivision that refines the expansion data for both $\mathscr R_D$ and $\mathscr R_T$, there is a commutative triangle
\[
\begin{tikzcd}
\mathscr R_D \arrow[rr,"\jmath"] \arrow[dr] &&
\mathscr R_T \arrow[dl]\\
& \mathscr U
\end{tikzcd}
\]
where $\mathscr U$ parametrizes prestable logarithmic maps to expansions of $\mathscr A$.
Strictness preserves the contact data, and the closed immersion
preserves stability. Then, the stack $\mathscr U$ is pure-dimensional
by \cite{Ran}*{Proposition 4.3.2} and by \cite{Ran}*{Section 4.5 and Theorem 4.5.1},
both maps to $\mathscr U$ carry the perfect
obstruction theories respectively, and their virtual classes are the
virtual pullbacks of $[\mathscr U]$.

Let $\pi_W$ be the map from universal curve $\mathcal C_W$ to $\mathscr R_D$, and let
$\varphi_W$ be the map from $\mathcal C_W$ to $W$ after contracting the target
expansion. Also, the contraction $\tau:\mathcal C_W\longrightarrow
\mathscr R_D\times_{\mathscr B_D}\mathcal C_D$ to the stabilized projected
curve satisfies 
\[
R\tau_*\mathcal O_{\mathcal C_W}
\simeq
\mathcal O_{\mathscr R_D\times_{\mathscr B_D}\mathcal C_D}.
\]
The projection formula and the derived base change give
\[
\bigl(R\pi_{W*}\varphi_W^*\rho^*N\bigr)^\vee[1]
\simeq
\mathbf Lp_D^*\mathbb E_\iota
\simeq
p_D^*\cV^\vee[1].
\]
Therefore, the normal sequence
\[
0\longrightarrow T_W(-\log(F_0+E_W))
\longrightarrow j^*T_{W_T}(-\log(F_{0,T}+E_T))
\longrightarrow\rho^*N
\longrightarrow0
\]
induces the upper distinguished
triangle in the following diagram
\[
\begin{tikzcd}[column sep=small,row sep=large]
\mathbf L\jmath^*\mathbb E_{\mathscr R_T/\mathscr U}
  \arrow[r]\arrow[d]
&
\mathbb E_{\mathscr R_D/\mathscr U}
  \arrow[r]\arrow[d]
&
p_D^*\cV^\vee[1]
  \arrow[r]\arrow[d]
&
(\mathbf L\jmath^*\mathbb E_{\mathscr R_T/\mathscr U})[1]
  \arrow[d]
\\
\mathbf L\jmath^*\mathbb L_{\mathscr R_T/\mathscr U}
  \arrow[r]
&
\mathbb L_{\mathscr R_D/\mathscr U}
  \arrow[r]
&
\mathbb L_{\mathscr R_D/\mathscr R_T}
  \arrow[r]
&
(\mathbf L\jmath^*\mathbb L_{\mathscr R_T/\mathscr U})[1]
\end{tikzcd}
\]
in which the lower row is the cotangent triangle. The vertical
morphisms come from the target cotangent triangle by
universal-map adjunction. Functoriality makes the diagram
commute including the connecting morphisms. 

Since $\cV$ is a vector bundle, the third relative obstruction theory
is perfect of amplitude $[-1,0]$. We therefore have a
compatible triple in the sense of
\cite{Manolache}*{Definition 4.5} and
\cite{Manolache}*{Corollary 4.9}, which gives
\[
\jmath^![\mathscr R_T]^\vir=[\mathscr R_D]^\vir
\]
where $\jmath^!$ is the base change of $\iota^!$.

Moreover, one has 
\[
\begin{tikzcd}[column sep=large,row sep=large]
\mathbf Lp_D^*\mathbb E_\iota
  \arrow[r,"\sim"]\arrow[d]
&
p_D^*\cV^\vee[1]
  \arrow[d]
\\
\mathbf Lp_D^*\mathbb L_{\mathscr B_D/\mathscr B_T}
  \arrow[r]
&
\mathbb L_{\mathscr R_D/\mathscr R_T}.
\end{tikzcd}
\]
The lower arrow is the morphism between the cotangent complexes for
\eqref{eq:vanishing-square}.

Since $T$ is convex, $\mathscr B_T$ is smooth of its
expected dimension. Lemma~\ref{lem:c1} shows
$(p_T)_*(\zeta_T\cap[\cR_T]^\vir)=0$ by a virtual-dimension counting argument.
Then, commutativity between virtual pullback and proper pushforward due to
\cite{Manolache}*{Theorem 4.1(i)} now gives
\[
\begin{aligned}
(p_D)_*(\zeta\cap[\mathscr R_D]^\vir)
&=(p_D)_*\jmath^!(\zeta_T\cap[\mathscr R_T]^\vir)\\
&=\iota^!(p_T)_*(\zeta_T\cap[\mathscr R_T]^\vir)\\
&=0.
\end{aligned}
\]
\end{proof}

In the following, we study the vanishing property of a $W$-vertex moduli space $\mathscr R_v$ after gluing in the degeneration formula (\ref{eq:relative-degeneration-explicit}) for any $v$ satisfying the assumptions of Lemma \ref{lem:vanishing}. Prior to that, we need to set up a bit of notation. 

For a $W$-vertex $v$ and a gluing edge $e$ linking $v$ to a $Q$-vertex, each map $f:C\rightarrow W$ (after contracting the expansion) maps the corresponding intersection point $x_e$ to $F_0$. Recall that each such a composition below is a constant map 
\[
C\xrightarrow{\,f\,}W\xrightarrow{\,g_W\,}\mathbb P^1,
\]
and thus the edge $e$ gives rise to an evaluation map 
\[
t_{v,e}:\mathscr R_v\longrightarrow\mathbb P^1,
\qquad
[f]\longmapsto g_W\bigl(f(x_e)\bigr).
\]
We denote the pullback \(t_{v,e}^*H\) of the point class from $\mathbb P^1$ by $H_{v, e}$.
Then, we can immediately observe the following lemma.

\begin{lemma}\label{prod:point-class}
    For a $W$-vertex $v$ and two edges $e, e'$ adjacent to $v$, the following identities hold
    \[
\begin{aligned}
H_{v,e}\,H_{v,e'}&=0,\\
\operatorname{ev}_0^*[X_{0,W}]\,H_{v,e}&=0 &&\text{if }v\text{ carries }p_0.
\end{aligned}
\]
\end{lemma}
\begin{proof}
    For the first equality, one has 
    \[
H_{v,e}H_{v,e'}
=t_v^*H\cdot t_v^*H
=t_v^*(H^2)=0.
\]

For the second equality, recall that $g_W(X_{0, W})=c_0$ is constant. Note that $t_v=g_W\circ\operatorname{ev}_0$ because $v$ carries $p_0$, and hence, 
\[
\operatorname{ev}_0^*[X_{0,W}]\,H_{v,e}
=\operatorname{ev}_0^*\!\left([X_{0,W}]\cdot g_W^*H\right).
\]
We can choose any distinct point $c\neq c_0$ to represent the class $H$. Therefore, $g_W^*H=[g_W^{-1}(c)]$ where $g_W^{-1}(c)$ does not intersect $X_{0,W}$. So, \([X_{0,W}]\cdot g_W^*H=0\), and we conclude that $\operatorname{ev}_0^*[X_{0,W}]\,H_{v,e}=0$.
\end{proof}

\begin{proposition}\label{prop:vanishing}
    For a rigid splitting graph $\Gamma$ containing a $W$-vertex $v$ satisfying Lemma \ref{lem:vanishing}, with the notation in (\ref{prop:relative-degeneration}), one has 
    \[
(u_\Gamma)_*\Delta_\Gamma^!(\operatorname{ev}_0^*[X_{0,\mathfrak V}]
\cap\boxtimes_w[\mathscr R_w]^{\mathrm{vir}}) =0,
\]
that is, such a graph contributes zero in the decomposition in Proposition \ref{prop:relative-degeneration} of $L_V$.
\end{proposition}
\begin{proof}
    First of all, in order to simplify the notation, we write $Z_\Gamma=\operatorname{ev}_0^*[X_{0,\mathfrak V}]
\cap\boxtimes_w[\mathscr R_w]^{\mathrm{vir}}$. Let $I$ be the edge set of $\Gamma$ and define $s:I\rightarrow V(\Gamma)$
write 
\(
M_\Gamma=\prod_{w\in V(\Gamma)}\mathscr R_w.
\)
Each edge $e=(w_1, w_2)\in I$ defines a map \(\operatorname{ev}_{F_0}:M_\Gamma\to(F_0\times F_0)^I\), and composing with the first projection map $\operatorname{pr}_D: F_0=D\times\mathbb P^1\rightarrow D$ gives the evaluation map $\operatorname{ev}_D:M_\Gamma\rightarrow (D\times D)^I$ at $D$. The gluing at $F_0$-level gives the glued space 
\[
\mathcal G_\Gamma
:=
M_\Gamma
\underset{(F_0\times F_0)^I}{\times}
F_0^I,
\]
which is defined in Section \ref{subsec:relative-degeneration}.
However, the space where $D$-evaluations agree is
\[
\mathcal G_{\Gamma,D}
:=
M_\Gamma
\underset{(D\times D)^I}{\times}
D^I.
\]
It is obvious that the projection $F_0\rightarrow D$ induces a closed immersion $i:\mathcal G_\Gamma\hookrightarrow\mathcal G_{\Gamma, D}$. Furthermore, note that we have a natural projection $p_D:\mathscr R_v\rightarrow\mathscr B_D$, thus if we set 
\[
M_\Gamma^{(v)}
:=
\mathscr B_D\times\prod_{w\ne v}\mathscr R_w,
\qquad
\mathcal G_{\Gamma,D}^{(v)}
:=
M_\Gamma^{(v)}
\underset{(D\times D)^I}{\times}
D^I,
\]
there is a map $q:\mathcal G_{\Gamma,D}\longrightarrow\mathcal G_{\Gamma,D}^{(v)}$.
With these definitions of maps, we have a factorization of $u_\Gamma$ as 
\[
\mathcal G_\Gamma
\xhookrightarrow{\ i\ }
\mathcal G_{\Gamma,D}
\xrightarrow{\ q\ }
\mathcal G_{\Gamma,D}^{(v)}
\xrightarrow{\ \bar u_\Gamma\ }
\mathscr C,
\qquad
u_\Gamma=\bar u_\Gamma\circ q\circ i.
\]
where $\bar u_\Gamma$ is the natural forgetful map. 

Notice that \([\Delta_{\mathbb P^1}]=H\otimes1+1\otimes H\), and thus for each choice of $s:I\rightarrow V(\Gamma)$ with $s(e)\in\partial e$ in which $\partial e$ represents the two endpoints, we have 
\[
\begin{aligned}
i_*\Delta_\Gamma^!Z_\Gamma
&=\Delta_{D,\Gamma}^!\left(
\prod_{e=(w,w')}(H_w+H_{w'})\cap Z_\Gamma\right)\\
&=\Delta_{D,\Gamma}^!\left(
  \prod_{e\in I}\left(\sum_{t\in\partial e}H_{t}\right)
  \cap Z_\Gamma\right)\\
&=\sum_{\substack{s:I\to V(\Gamma)\\s(e)\in\partial e\ \forall e\in I}}
  \Delta_{D,\Gamma}^!\left(
  \bigl(\zeta_s\cap[\mathscr R_v]^{\mathrm{vir}}\bigr)
  \boxtimes\alpha_s\right)
\end{aligned}
\]
where $\Delta_{D,\Gamma}:D^I\longrightarrow(D\times D)^I$ and 
\[
\zeta_{t,s}
:=
\left(\prod_{\substack{e\in I\\s(e)=t}}H_{t,e}\right)
\begin{cases}
\operatorname{ev}_0^*[X_{0,\mathfrak V}],&p_0\text{ lies on }t,\\
1,&\text{otherwise}.
\end{cases}
\]
Here, $\zeta_s:=\zeta_{v,s}$ with $v$ we fix at the very beginning and $\alpha_s:=
\mathop{\boxtimes}_{t\ne v}
\bigl(\zeta_{t,s}\cap[\mathscr R_t]^{\mathrm{vir}}\bigr)$.

Then, by the preceding lemma, each $\zeta_s$ is either $0$ or one of \(1\), \(H_{v,e}\), and \(\operatorname{ev}_0^*[X_{0,W}]\). So, by Lemma \ref{lem:vanishing}, we have 
\[
\begin{aligned}
&(u_\Gamma)_*\Delta_\Gamma^! Z_\Gamma
=(\bar u_\Gamma)_*q_*i_*\Delta_\Gamma^! Z_\Gamma\\[3pt]
&\quad=
\sum_{\substack{s:I\to V(\Gamma)\\ s(e)\in\partial e\ \forall e\in I}}
(\bar u_\Gamma)_*q_*\Delta_{D,\Gamma}^!
\left(
  \bigl(\zeta_s\cap[\mathscr R_v]^{\mathrm{vir}}\bigr)
  \boxtimes\alpha_s
\right)\\[3pt]
&\quad=
\sum_{\substack{s:I\to V(\Gamma)\\ s(e)\in\partial e\ \forall e\in I}}
(\bar u_\Gamma)_*\Delta_{D,\Gamma}^!
\left(
  \underbrace{(p_D)_*\bigl(\zeta_s\cap[\mathscr R_v]^{\mathrm{vir}}\bigr)}_{=\,0}
  \boxtimes\alpha_s
\right)
=0.
\end{aligned}
\]
This completes the proof of our vanishing claim. 
\end{proof}

\begin{remark}\label{rem:stable-ghosts}
Proposition \ref{prop:vanishing} handles the case where $(\delta_D, q)$ marks stable logarithmic maps to $D$, that is, either $\delta_D=0$ with $q\geq 3$ or $\delta_D\neq0$. In either case, Proposition \ref{prop:vanishing} tells us that the corresponding glued virtual cycles will be $0$ after pushing forward to $\mathscr C$, so they will not contribute to the decomposition of $L_V$. Hence, the only case left to compute is 
\[
\delta_D=0,\qquad q\leq2,
\]
which will be discussed in the following sections. 
\end{remark}

\section{The contributions in a fiber}\label{sec:caps}
Now, as Remark \ref{rem:stable-ghosts} discusses, we now focus on the situation that $\delta_D=0$ and $q\leq 2$, that is, we only consider those stable logarithmic maps which carry at most two markings including gluing markings and become constant after projection to $D$. 


These glued maps have images which lie in a fiber surface $S_d$ of
$\rho:W\to D$ which is just the blowup of $\mathbb P^2$ at two distinct points $p_1, p_2$, i.e. $S_d\simeq \operatorname{Bl}_{p_1, p_2}\mathbb P^2$. 
Then let $s$ be the ($-1$)-curve that is contained in $X_{0,W}$ as well, that is, $s=X_{0,W}\cap S_d$ which has a normal line bundle $\mathcal O_{\mathbb P^1}(-1)$. Let $\ell_W$
be the strict transform of the line through $p_1, p_2$, which is another ($-1$)-curve of $S_d$.

If $a$ is the total contact order with $F_0$ and $k$ the total
contact order with $E_W$, the curve class in question is of the form
\begin{equation}\label{eq:fiber-class}
a s+k\ell_W.
\end{equation}
In the following proposition, a coefficient represents the multiplicity of the virtual pushforward to $D$.

\begin{proposition}[Fiber coefficients]\label{prop:caps}
The only nonzero $W$-vertex contributions are as follows.
\begin{enumerate}
\item A vertex with one contact marking $p_i$ ($1\leq i\leq n$) having contact order $k_i$ with $E_W$ and one gluing marking of order $a$ where $a=k_i$, which contributes a coefficient $1/k_i$ after pushing forward to $D$.
\[
\begin{tikzcd}[column sep=3em, cells={nodes={inner sep=0pt}}]
  \overset{Q}{\circ}
    \arrow[r,dash,teal,thick,"{a=k_i}",""'{pos=.93,text=blue!65!black}] &
  \overset{W}{\bullet}
    \arrow[r,dash,dashed,purple] &
  \textcolor{purple}{p_i\text{ on }E_W,\ \text{order }k_i}
\end{tikzcd}
\]
\item A vertex with one gluing marking of order $a$ and no original
markings contributes a coefficient
$(-1)^{a+1}/a^2$.
\[
\begin{tikzcd}[column sep=5em, cells={nodes={inner sep=0pt}}]
  \overset{Q}{\circ} \arrow[r,dash,teal,thick,"a"] &
  \overset{W}{\bullet}
\end{tikzcd}
\]
\item The vertex in (2) with the additional marking $p_0$
and insertion $\ev_0^*[X_{0,W}]$ contributes a coefficient
$(-1)^a/a$.
\[
\begin{tikzcd}[column sep=3em, cells={nodes={inner sep=0pt}}]
  \overset{Q}{\circ} \arrow[r,dash,teal,thick,"a"] &
  \overset{W}{\bullet} \arrow[r,dash,dashed,purple] &
  \textcolor{purple}{p_0\;\bigl(\operatorname{ev}_0^*[X_{0,W}]\bigr)}
\end{tikzcd}
\]
\end{enumerate}
An unmarked vertex with two gluing edges, or a one-edge vertex
with an additional ordinary interior marking other than $p_0$,
has zero pushforward.
\end{proposition}
\begin{proof}
For (1), fix a $W$-vertex $v$ with gluing marking $x$ and an edge $e=(x,x')$ linking to another $Q$-vertex with $x'$. Since $F_0=D\times\mathbb P^1$, there is a pair of evaluations at $x$
\[
\operatorname{ev}_x=(\epsilon_x,t_x)\colon
\mathscr R_v\longrightarrow D\times\mathbb P^1.
\]
Write $t_{x'}:\mathscr R_v\rightarrow\mathbb P^1$ for the evaluation at $x'$. The gluing gives the class 
\[
(t_{x'},t_x)^*[\Delta_{\mathbb P^1}]
=t_{x'}^*H+t_x^*H,
\qquad H=[c]\in A^1(\mathbb P^1),
\]
where $c$ is the point class in $\mathbb P^1$. Therefore, for this $W$-vertex $v$, the gluing will produce terms 
\[
[\mathscr R_v]^{\mathrm{vir}}
\qquad\text{and}\qquad
t_x^*H\cap[\mathscr R_v]^{\mathrm{vir}}
\]
in the degeneration formula. 
Notice that $\operatorname{vdim}\mathscr R_v=\operatorname{dim}D + 1=\operatorname{dim}F_0$, we have 
\[
(\operatorname{ev}_x)_*[\mathscr R_v]^{\mathrm{vir}}
=b_{a,k_i}[F_0].
\]
To compute the coefficient $b_{a,k_i}$, let us fix $(d,c)\in F_0=D\times\mathbb P^1$. As we have seen earlier that both projection to $D$ and the pencil map to $\mathbb P^1$ are constant the connected source curve, so the image of $\operatorname{ev}_x$ lies in $\rho^{-1}(d)\cap g_W^{-1}(c)\simeq\mathbb P^1$. The divisors $F_0$ and $E_W$ of $W$ cut out two points on this curve and we denote them by $0$ and $\infty$ respectively. This enforces $a=k_i$ since $a$ is the contact order with $0$ and $k_i$ is the contact order with $\infty$ but they both equal the degree of a branched cover $C\rightarrow\mathbb P^1$, thus the cover is given by $z\mapsto z^{k_i}$ and obviously unobstructed. Therefore, $b_{a, k_i}=1/k_i$.
Finally, note that \(\epsilon_x=\operatorname{pr}_D\circ\operatorname{ev}_x\), we have 
\[
\begin{aligned}
(\epsilon_x)_*[\mathscr R_v]^{\mathrm{vir}}
&=b_{a,k_i}(\operatorname{pr}_D)_*[D\times\mathbb P^1]
=0,\\[3pt]
(\epsilon_x)_*\!\left(t_x^*H\cap[\mathscr R_v]^{\mathrm{vir}}\right)
&=b_{a,k_i}(\operatorname{pr}_D)_*[D\times\{c\}]
=b_{a,k_i}[D]=1/k_i[D].
\end{aligned}
\]

For (2), still fix a $W$-vertex $v$ and let $x$ be the gluing marking of contact order
$a>0$ along $F_0$, and we write $\ev_x=(\epsilon_x,t_x)$ as defined in (1).
Since it is assumed that there is no additional markings at $v$, the curve class in
$S_d=\rho^{-1}(d)$ is $a s$ where recall that $s=X_{0,W}\cap S_d$. Let $f:C\to S_d$ be a map obtained
by contracting the target expansion. Under the blow-down
$\nu:S_d\to\PP^2$, the curve $s$ is the exceptional curve over some point $p$.
Hence,
\[
\deg(\nu\circ f)^*\OO_{\PP^2}(1)
=\nu^*c_1(\OO_{\PP^2}(1))\cdot (a s)=0.
\]
Note that $\deg(\nu\circ f)^*\OO_{\PP^2}(1)$ is nonnegative on each component of $C$, so $\nu\circ f$
is constant on the connected source curve $C$. Since $a>0$, its value must
be one of the two blowup centers $p,q$. The exceptional curves
over $p$ and $q$ have independent classes, and $f_*[C]=a s$.
It follows that $f(C)\subset s$. 

Since $s$ meets $F_0$ in exact one point and is disjoint from $E_W$,
every source component getting fully mapped into an $F_0$-rubber component maps into the $\mathbb P^1$ fiber over the point $s\cap F_0$, and no any
$E_W$-level expansions occur due to the stability issue. Let $Y_\infty$ be the infinity
section of $\PP_Y(\OO_Y(-D_\infty)\oplus\OO_Y)$ and let
$U=W\setminus b_W^{-1}(Y_\infty)$. Then the exceptional curve $s$ lies in $U$,
so these maps take values in expansions of the pair
\[
(U,F_0\cap U)
\simeq
\bigl(\operatorname{Tot}(\OO_Y(-D_\infty)),\,
      D_0\times\AA^1\bigr).
\]
Also, we have shown that for any stable logarithmic map with expansions \(\widetilde f:C\to S_d^{\mathrm{exp}}\) with \(f:C\to S_d\) the corresponding map obtained by composing with the contraction \(S_d^{\mathrm{exp}}\to S_d\), one has $f(C)\subset s$. Therefore, the open pair in $S_d$ containing $s$ is as follows
\[
(U\cap S_d,F_0\cap U\cap S_d)
\simeq
\bigl(\operatorname{Tot}(\mathcal O_{\mathbb P^1}(-1)),
      \mathcal O_{\mathbb P^1}(-1)|_q\bigr),
\]where $q=F_0\cap s$.
The fiber calculation in the proof of \cite{vGGR}*{Proposition 2.4} gives 
\[
(\epsilon_x)_*[\mathscr R_v]^\vir
=\frac{(-1)^{a+1}}{a^2}[D].
\]

For (3), this case is just case (2) with one additional divisor insertion $\operatorname{ev}_0^*[X_{0,W}]$ at $p_0$.
Note that $i:Y\simeq X_{0,W}\hookrightarrow W$ gives rise to
\[
i^*[X_{0,W}]
=c_1\!\left(\mathcal O_W(X_{0,W})|_Y\right)
=c_1\!\left(\mathcal O_Y(-D_\infty)\right)
=-[D_\infty],
\]
which gives the self-intersection of $X_{0,W}$.
Hence,
\[
X_{0,W}\cdot s
=\deg\!\left(\mathcal O_W(X_{0,W})|_s\right)
=\deg\!\left(\mathcal O_Y(-D_\infty)|_s\right)
=-1.
\]
Therefore, one has 
\[
\begin{aligned}
(\epsilon_x)_*\!\left(
\operatorname{ev}_{p_0}^*[X_{0,W}]
\cap[\mathscr R_v^+]^{\mathrm{vir}}\right)
&=(X_{0,W}\cdot a s)\,
  (\epsilon_x)_*[\mathscr R_v]^{\mathrm{vir}}\\
&=(-a)\frac{(-1)^{a+1}}{a^2}[D]\\
&=\frac{(-1)^a}{a}[D].
\end{aligned}
\]
\end{proof}

\begin{remark}[The other vertex types]\label{rem:exhaustion}
If a bipartite graph contains vertices over the $Q$-side, each $W$-vertex has at least one gluing marking thanks to Proposition \ref{prop:caps}. Lemma \ref{lem:vanishing} eliminates possibilities that the projected curve class in $D$ is non-constant and it has at least three markings in total. Then, Proposition \ref{prop:caps} shows that only vertices having one gluing marking can contribute non-trivially. Also, the contributing graph cannot fully lie over the $W$-side. Therefore, we can now claim that a graph contributing non-trivially must be star-shaped. 
\end{remark}

\section{The degeneration graph and the conclusion}
\label{sec:assembly}

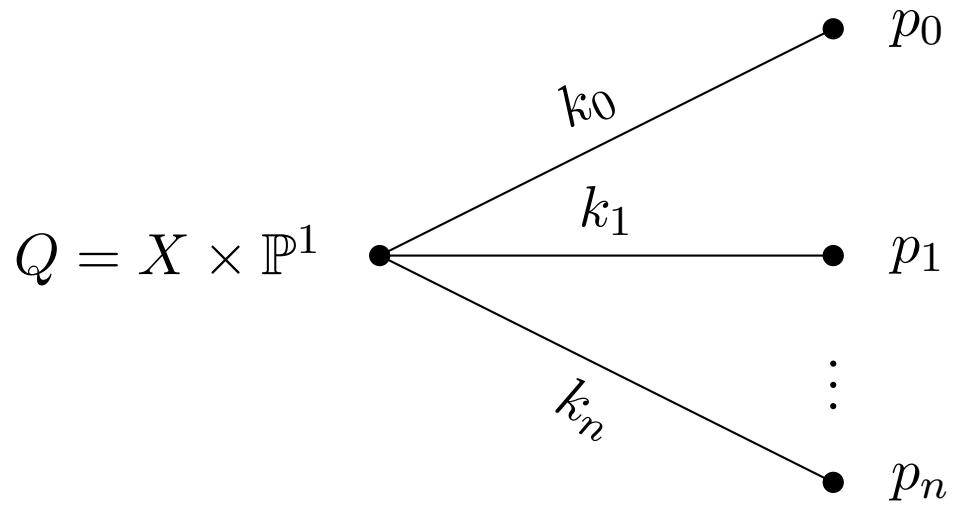
\begin{figure}[ht]
\centering
\begin{tikzpicture}[x=1cm,y=1cm]
\fill (0,0) circle (2pt);
\node[left=7pt] at (0,0) {$Q=X\times\PP^1$};
\draw (0,0)--(3,1.5) node[midway,above,sloped] {$k_0$};
\draw (0,0)--(3,0) node[midway,above] {$k_1$};
\draw (0,0)--(3,-1.5) node[midway,below,sloped] {$k_n$};
\fill (3,1.5) circle (2pt);
\fill (3,0) circle (2pt);
\fill (3,-1.5) circle (2pt);
\node[right=7pt] at (3,1.5) {$p_0$};
\node[right=7pt] at (3,0) {$p_1$};
\node[right=7pt] at (3,-1.5) {$p_n$};
\node at (3,-0.75) {$\vdots$};
\end{tikzpicture}
\caption{The only contributing graph type, drawn for $n\geq2$.
The vertices on the right side map to $W$ and carry the indicated
markings. All additional interior markings lie on $Q$.}
\label{fig:star}
\end{figure}

Now, we are in a position to prove Proposition~\ref{prop:relative}, which is the comparison
between the moduli stacks of stable logarithmic maps to expansions for the two distinct smooth pairs. Lemma \ref{rem:exhaustion} already implies that the surviving terms in the degeneration formula of $L_V$ must be given by gluing the vertex moduli spaces on star-shaped graphs. We will now show
that only the particular labeled star-shaped graph up to isomorphisms depicted in Figure~\ref{fig:star} contributes non-trivially.
One crucial point to keep in mind is that a vertex over $Q=X\times\PP^1$
has one free $\PP^1$ parameter. Hence, integrating over this component
requires exactly one point insertion to result in a non-zero outcome.

\begin{proposition}[The contributing graph]\label{prop:graph}
The contribution of a degeneration graph vanishes unless
it has one vertex over $Q$ and $n+1$ leaves over $W$, each of which is joined to $Q$ by one edge.
For each $0\leq i\leq n$, the leaf carrying $p_i$ is linked to the $Q$-vertex by an edge of contact order $k_i$
All additional interior markings lie on $Q$.
\end{proposition}

\begin{proof}
Let $\Gamma$ be a graph with nonzero contribution. By
Proposition~\ref{prop:caps} and Remark~\ref{rem:exhaustion},
$\Gamma$ is star-shaped with a central vertex $v_Q$ over $Q$. Each $p_i$, for
$1\leq i\leq n$, has positive contact with $E_W$ by our fixed contact order profile, so it must be attached to a
$W$-leaf as a marking. Proposition~\ref{prop:caps}(1) shows that these leaves are
distinct and that the edge $e_i$ incident to the leaf carrying $p_i$
has order $a_{e_i}=k_i$. Therefore, a graph contributing to the degeneration formula has at least $n$ distinct $W$-leaves.

Again, by Lemma \ref{lem:vanishing}, the projected curve class is zero in $D$, and thus $v_Q$ carries the class whose pushforward becomes our fixed curve class $\beta$ in $X$.
Let $\mathscr R_X$ be the moduli space of stable logarithmic maps with expansions to $(X,D)$ with
the data determined at $v_Q$. Note that the expansions of $(Q,F_0)$ are the products of
expansions of $(X,D)$ with $\PP^1$. More importantly, Lemma \ref{lem:pencil} says 
\[
g^*H\cdot\lambda
=X_\infty\cdot\gamma-E\cdot\lambda
=K-K=0.
\]
Here \(X_\infty\cdot\gamma=K\) because \(\gamma=\beta+Kf\).
Therefore, each map to $Q$ in the moduli stack $\mathscr R_{v_Q}$ at $v_Q$ is a constant $c$ along $\PP^1$. This $\PP^1$-factor is unobstructed
since $H^1(C,\OO_C)=0$. Hence,
\[
\mathscr R_{v_Q}\simeq\mathscr R_X\times\PP^1,
\qquad
[\mathscr R_{v_Q}]^\vir
=[\mathscr R_X]^\vir\boxtimes[\PP^1].
\]
Let us order the factors $\PP^1$ of each gluing diagonal as $Q,W$ so that
\[
[\Delta_{\PP^1}]=H\otimes 1+1\otimes H,
\]
that is, \(H\otimes 1\) places the point class on \(Q\), and \(1\otimes H\) places it on \(W\).
For a leaf carrying $p_i$, $i>0$, the proof of Proposition \ref{prop:caps} implies that the virtual cycle at this leaf pushes forward to zero unless the leaf acquires a point class insertion $H$ at the gluing marking, hence only the summand \(1\otimes H\) contributes. 
For a leaf with no $E_W$-contact marking, the parametrized maps have constant value $c_0$ along $\PP^1$,
so only $H\otimes 1$ survives. Let $r$ be the number of leaves
with no original marking. Each of these $r$ leaves therefore
supplies one $H$ on $Q$. However, for the marking $p_0$, it supplies one more factor of $H$. 
If it is attached at $v_Q$, then $[X_{0,Q}]=\pr_2^*H$, and if it lies on
$W$-side, the edge actually supplies $H\otimes 1$ which survives, whereas the summand in the degeneration formula containing $1\otimes H$ vanishes.

Now, we see that the resulting point class on $Q$ is $H^{r+1}$. However, 
\[
H^{r+1}=
\begin{cases}
H,&r=0,\\
0,&r\ge 1.
\end{cases}
\]
Writing $\pr_1:\mathscr R_X\times\PP^1\to\mathscr R_X$ for the
projection, we therefore obtain
\[
\begin{aligned}
(\pr_1)_*\!\left(
  \pr_2^*H^{r+1}\cap[\mathscr R_{v_Q}]^\vir
\right)
&=\left(\int_{\PP^1}H^{r+1}\right)[\mathscr R_X]^\vir\\
&=
\begin{cases}
[\mathscr R_X]^\vir,& r=0,\\
0,& r\geq 1.
\end{cases}
\end{aligned}
\]
Therefore, every contributing graph has $r=0$.

It remains to locate $p_0$. Let $\delta_Q$ be the curve class carried at
$v_Q$. All $W$-leaves give out stable logarithmic maps that project to points in $D$, so
$(\pr_X)_*\delta_Q=\beta$. Hence, we have 
\[
\sum_{e\ni v_Q}a_e
=F_0\cdot\delta_Q
=D\cdot\beta
=d.
\]
If $p_0$ were laid on $Q$, the only edges would be $e_1,\ldots,e_n$,
whose orders sum to $\sum_{i=1}^n k_i=d-k_0<d$ since we assume that $k_0> 0$. Hence, $p_0$
lies on a $W$-leaf. Its edge $e_0$ has order
\[
a_{e_0}
=d-\sum_{i=1}^n a_{e_i}
=d-\sum_{i=1}^n k_i
=k_0.
\]
The leaves are therefore indexed by $p_0,\ldots,p_n$. Finally, again by
Proposition~\ref{prop:caps}, all additional interior markings
are attached at $v_Q$.
\end{proof}

\begin{proof}[Proof of Proposition~\ref{prop:relative}, for connected $D$]
Apply the degeneration formula stated in Section~\ref{sec:degeneration} to $L_V$.
Lemma~\ref{lem:vanishing} and Propositions~\ref{prop:caps}
and~\ref{prop:graph} show that only the graph above non-trivially contributes.

Projection to $D\in X$ contracts every $W$-leaf and transfers its
corresponding marking $p_i$ for some $i$ to the relative marking to the divisor $D$ in ($X,D$). All retained $D$-evaluations are preserved. Stable logarithmic maps to $(Q=X\times\PP^1,D\times\PP^1)\cup (W, F_0+E_W)$ marked by this graph will have the contact order $k_i$ at each $p_i$ after they we project these maps to $(X,D)$. So, obviously, 
pushing the virtual class $[\mathscr R_X]^{\operatorname{vir}}$ forward via the map 
$\pi_{\operatorname{rel}}:\mathscr R_X\rightarrow \mathscr C$
gives exactly
$R_{\mathrm{rel}}$.

For the vertex carrying $p_i$, $i>0$, the coefficient $1/k_i$
cancels the edge multiplicity $k_i$. For the vertex carrying
$p_0$, the coefficient is $(-1)^{k_0}/k_0$ and the edge
multiplicity is $k_0$. Their product is
\[
\left(\prod_{i=1}^n k_i\frac1{k_i}\right)
\left(k_0\frac{(-1)^{k_0}}{k_0}\right)=(-1)^{k_0}.
\]
The remaining integral is $\int_{\PP^1}H=1$ on the $Q$-vertex.
The original marking labels distinguish all the leaves, so
the automorphism group of the graph preserving all markings with labels is trivial.
It follows that
\[
L_V=(-1)^{k_0}R_{\mathrm{rel}}.
\]
\end{proof}

\begin{proof}[Proof of Theorem~\ref{thm:main}, for connected $D$]
Proposition~\ref{prop:smooth} gives $L=L_V$ and
Proposition~\ref{prop:relative} gives
$L_V=(-1)^{k_0}R_{\mathrm{rel}}$.
A smooth log-relative comparison for $(X,D)$ gives
$R_{\mathrm{rel}}=R$. Thus $L=(-1)^{k_0}R$, which is
equivalent to \eqref{eq:main}. The same identities for
disconnected $D$ are established in Section~\ref{sec:endpoints}.
\end{proof}

\section{The maximal-contact case and disconnected divisors}
\label{sec:endpoints}

\subsection{The maximal-contact endpoint}
If $K=0$, then $n=0$ because every contact order $k_i$ with $1\leq i\leq n$ is positive. In this case, the graph from Proposition \ref{prop:graph} literally becomes the only contributing graph exhibited in \cite{vGGR} in their setting, which is a line segment with two vertices where one is $Q$-vertex and the other is $W$-vertex. This is literally the maximal contact case and the curve class is literally $\gamma=(i_0)_*\beta$ where $i_0:X\hookrightarrow P$ is the zero section.
The $W$-vertex carries $p_0$ in this case and it has degree $d=k_0$,
and gives $(-1)^d$ after including the edge multiplicity.

If we further assume that $D$ is a nef divisor in $X$. Then $X_\infty$ is nef on $P$, and thus we have 
\[
X_\infty\cdot\gamma=0,\qquad X_0\cdot\gamma=-d<0.
\]
The negativity of intersection with $X_0$ together with the connectedness of the source curves forces the whole map to stay in $P\setminus X_\infty=\operatorname{Tot}(\mathcal O_X(-D))$. Furthermore, our comparison contains the insertion $\operatorname{ev}_0^*[X_0]$. Since $i_0^*[X_0]=-[D]$, and forgetting the marking $p_0$ gives a factor $-d$ based on the divisor equation. Therefore, our comparison in Equality (\ref{eq:main}) in Theorem \ref{thm:main} will come down to the main result in \cite{vGGR}
\[
F_*[\overline{\mathcal M}_{0,(d)}(X(\log D),\beta)]^{\mathrm{vir}}
=
(-1)^{d+1}d\,
[\overline{\mathcal M}_{0,0}(\operatorname{Tot}(\mathcal O_X(-D)),\beta)]^{\mathrm{vir}}.
\]
Therefore, our Theorem \ref{thm:main} recovers the result of \cite{vGGR}.

For $n>0$, the problem we consider here has several positive relative
markings, and it is not the maximal-contact problem for $(X,D)$. Our 
Theorem \ref{thm:main} removes one of these contacts and retains the
others through \eqref{eq:contacts}.

\subsection{A disconnected smooth divisor}
Suppose $D=\coprod_\alpha D_\alpha$. Its components are
disjoint.
Fix a $D$-component assignment for each marking $p_i$ for $1\leq i\leq n$, and
let $K_\alpha$ be the sum of the contact orders
assigned to $D_\alpha$, thus $\sum_\alpha K_\alpha=K$.
Run the same construction as we did earlier for constructing $(V,E)$, and let $E_\alpha$ be the
exceptional divisor over $D_\alpha$ and $e_\alpha$ the fiber
class of $E_\alpha\to D_\alpha$.
Using the same lifted zero-section identification for $\beta$
and the same general-fiber class $f_V$, the lifted class is
\[
\lambda=\beta+Kf_V-\sum_\alpha K_\alpha e_\alpha.
\]
The degeneration gives rise to a central fiber comprised of one component $Q=X\times\PP^1$ and a
component $W_\alpha$ for each $D_\alpha$.

Apply Lemma~\ref{lem:vanishing} to each $W_\alpha$, then the corresponding proof applies componentwise. Note that we have
\[
D_\alpha\cdot\beta-K_\alpha=
\begin{cases}
k_0,&\alpha=\alpha_0,\\
0,&\alpha\ne\alpha_0,
\end{cases}
\]
where $D_{\alpha_0}$ is assumed to carry $p_0$.
These are precisely the component assignments of the original
relative problem, including the component of $p_0$.
Summing the compatible assignments proves
Proposition~\ref{prop:relative} for disconnected $D$.
Therefore, Theorem \ref{thm:main} holds for a disconnected $D$ as well. 

\begin{remark}
Note that Theorem \ref{thm:main} and its proof concern a smooth pair with no nefness condition imposed and the auxiliary
SNC pair with diagonal contacts \eqref{eq:contacts}.
It does not assert a correspondence for arbitrary SNC pairs or a higher-genus comparison.
\end{remark}

\bibliographystyle{amsxport}
\bibliography{main}

\clearpage
\begin{bibdiv}
\begin{biblist}

\bib{ACFW}{article}{
  author={Abramovich, Dan},
  author={Cadman, Charles},
  author={Fantechi, Barbara},
  author={Wise, Jonathan},
  title={Expanded degenerations and pairs},
  journal={Comm. Algebra},
  volume={41},
  date={2013},
  pages={2346--2386},
  eprint={arXiv:1110.2976}
}

\bib{AF}{article}{
  author={Abramovich, Dan},
  author={Fantechi, Barbara},
  title={Orbifold techniques in degeneration formulas},
  journal={Ann. Sc. Norm. Super. Pisa Cl. Sci. (5)},
  volume={16},
  date={2016},
  pages={519--579},
  eprint={arXiv:1103.5132}
}

\bib{AMW}{article}{
  author={Abramovich, Dan},
  author={Marcus, Steffen},
  author={Wise, Jonathan},
  title={Comparison theorems for Gromov--Witten invariants of smooth
    pairs and of degenerations},
  journal={Ann. Inst. Fourier},
  volume={64},
  date={2014},
  pages={1611--1667},
  eprint={arXiv:1207.2085}
}

\bib{AW}{article}{
  author={Abramovich, Dan},
  author={Wise, Jonathan},
  title={Birational invariance in logarithmic Gromov--Witten theory},
  journal={Compos. Math.},
  volume={154},
  date={2018},
  pages={595--620},
  eprint={arXiv:1306.1222}
}

\bib{BNR}{article}{
  author = {Battistella, Luca},
  author = {Nabijou, Navid},
  author = {Ranganathan, Dhruv},
  title  = {Gromov-Witten theory via roots and logarithms},
  eprint = {arXiv:2203.17224}

}

\bib{Fulton}{book}{
  author={Fulton, William},
  title={Intersection Theory},
  edition={2},
  publisher={Springer},
  date={1998}
}

\bib{Li}{article}{
  author={Li, Jun},
  title={A degeneration formula of GW-invariants},
  journal={J. Differential Geom.},
  volume={60},
  date={2002},
  pages={199--293},
  eprint={arXiv:math/0110113}
}

\bib{Man}{article}{
  author={Manolache, Cristina},
  title={Virtual pull-backs},
  journal={J. Algebraic Geom.},
  volume={21},
  date={2012},
  pages={201--245},
  eprint={arXiv:0805.2065}
}

\bib{vGGR}{article}{
  author={van Garrel, Michel},
  author={Graber, Tom},
  author={Ruddat, Helge},
  title={Local Gromov--Witten invariants are log invariants},
  journal={Adv. Math.},
  volume={350},
  date={2019},
  pages={860--876},
  eprint={arXiv:1712.05210}
}

\bib{Ran}{article}{
  author={Ranganathan, Dhruv},
  title={Logarithmic Gromov--Witten theory with expansions},
  journal={Algebraic Geometry},
  volume={9},
  date={2022},
  pages={714--761},
  eprint={arXiv:1903.09006}
}


\bib{TY20b}{article}{
    author = {Tseng, Hsian-Hua},
    author = {You, Fenglong},
     TITLE = {A mirror theorem for multi-root stacks and applications},
   JOURNAL = {Selecta Math. (N.S.)},
    VOLUME = {29},
      YEAR = {2023},
    NUMBER = {1},
     PAGES = {Paper No. 6, 33},
      ISSN = {1022-1824,1420-9020},
       DOI = {10.1007/s00029-022-00809-8},
       URL = {https://doi.org/10.1007/s00029-022-00809-8},
}


\bib{TY20c}{article}{
    AUTHOR = {Tseng, Hsian-Hua}, 
    author = {You, Fenglong},
     TITLE = {A {G}romov-{W}itten theory for simple normal-crossing pairs
              without log geometry},
   JOURNAL = {Comm. Math. Phys.},
    VOLUME = {401},
      YEAR = {2023},
    NUMBER = {1},
     PAGES = {803--839},
      ISSN = {0010-3616,1432-0916},
       DOI = {10.1007/s00220-023-04656-2},
       URL = {https://doi.org/10.1007/s00220-023-04656-2},
}


\bib{WY}{article}{
  author={Wang, Yu},
  author={You, Fenglong},
  title={Gromov--Witten theory beyond maximal contacts},
  date={2024},
  eprint={arXiv:2403.17200},
  note={Theorem 1.3 in the supplied manuscript}
}

\bib{Chen}{article}{
  author={Chen, Qile},
  title={Stable logarithmic maps to Deligne--Faltings pairs I},
  journal={Ann. of Math. (2)},
  volume={180},
  date={2014},
  number={2},
  pages={455--521},
  eprint={arXiv:1008.3090}
}

\bib{AC}{article}{
  author={Abramovich, Dan},
  author={Chen, Qile},
  title={Stable logarithmic maps to Deligne--Faltings pairs II},
  journal={Asian J. Math.},
  volume={18},
  date={2014},
  number={3},
  pages={465--488},
  eprint={arXiv:1102.4531}
}

\bib{GS13}{article}{
  author={Gross, Mark},
  author={Siebert, Bernd},
  title={Logarithmic Gromov--Witten invariants},
  journal={J. Amer. Math. Soc.},
  volume={26},
  date={2013},
  number={2},
  pages={451--510},
  eprint={arXiv:1102.4322}
}

\bib{KLR}{article}{
 author = {Kim, Bumsig},
 author = {Lho, Hyenho},
 author = {Ruddat, Helge},
 title = {The degeneration formula for stable log maps},
 eprint = {arXiv:1803.04210}

}

\end{biblist}
\end{bibdiv}

\if{

}\fi

\end{document}